%% file: article.tex
\documentclass[onefignum,onetabnum,oneeqnum]{siamart171218}

\input{shared}

\ifpdf
\hypersetup{
  pdftitle={}Vlasov-Poisson system tackled by particle simulation utilising boundary element methods,
  pdfauthor={Torsten Kessler, Sergej Rjasanow and Steffen Weisser}
}
\fi

\begin{document}

\maketitle

\begin{abstract}
This paper presents a grid-free simulation algorithm for the fully
three-dimensional Vlasov--Poisson system for collisionless
electron plasmas. We employ a standard
particle method for the numerical approximation of the distribution function.
Whereas the advection of the particles is grid-free by its very nature,
the computation of the acceleration involves the solution of the non-local
Poisson equation. To circumvent a volume mesh, we utilise the
Fast Boundary Element Method, which reduces the three-dimensional
Poisson equation to a system of linear equations on its two-dimensional boundary.
This gives rise to fully populated matrices which are approximated
by the $\mathcal H^2$-technique,
reducing the computational time from quadratic to linear complexity.
The approximation scheme based on interpolation has shown to be robust and
flexible, allowing a straightforward generalisation to vector-valued
functions.
In particular, the Coulomb forces acting on the particles are computed
in linear complexity.
In first numerical tests, we validate our approach with the help of classical
non-linear plasma phenomena. Furthermore, we show that our method
is able to simulate electron plasmas in complex three-dimensional domains with mixed
boundary conditions in linear complexity.
\end{abstract}

\begin{keywords}
Vlasov--Poisson system, simulation of plasmas, particle method,
Boundary Element Method, hierarchical approximation
\end{keywords}

\begin{AMS}
35Q83, 65Z05, 65N75, 65N38, 68W25
\end{AMS}

\input{introduction}
\input{vlasov-system}
\input{bem}
\input{hierarchical}
\input{implementation}
\input{numerics}
\bibliographystyle{siamplain}
\bibliography{literature}
\end{document}

%% file: shared.tex
\usepackage{lipsum}
\usepackage{subcaption}
\usepackage{amssymb}
\usepackage{amsfonts}
\usepackage{graphicx}
\usepackage{epstopdf}
\usepackage{algorithmic}
\ifpdf
  \DeclareGraphicsExtensions{.eps,.pdf,.png,.jpg}
\else
  \DeclareGraphicsExtensions{.eps}
\fi
\usepackage{bbm}

\newsiamremark{remark}{Remark}
\newsiamremark{hypothesis}{Hypothesis}
\crefname{hypothesis}{Hypothesis}{Hypotheses}
\newsiamthm{claim}{Claim}

\headers{Vlasov-Poisson system tackled by BEM}{T. Ke\ss ler, S. Rjasanow, and S. Wei\ss er}

\title{Vlasov-Poisson system tackled by particle simulation utilising Boundary Element Methods\thanks{Submitted to the editors DATE.}}

\author{Torsten Ke\ss ler\thanks{Department of Mathematics, Saarland University, 66041 Saarbr\"ucken, Germany (\email{kessler@num.uni-sb.de}).}
\and Sergej Rjasanow\thanks{Department of Mathematics, Saarland University, 66041 Saarbr\"ucken, Germany (\email{rjasanow@num.uni-sb.de}).}
\and Steffen Wei\ss er\thanks{Department of Mathematics, Saarland University, 66041 Saarbr\"ucken, Germany (\email{weisser@num.uni-sb.de}).}
}

\usepackage{amsopn}

\DeclareMathOperator{\Span}{span}


%% file: introduction.tex
\section{Introduction}\label{sec:introduction}

The rapid increase of computational power
in the last years due to massively parallel machines
like clusters or GPUs has opened up the possibility of 
handling complex problems for a broad range of applications
utilising classical particle methods.
Readily implemented in a computer program, they
are extensible and applicable to computational problems in biology, chemistry
and physics.

Particle methods for the simulation of collisionless plasmas has been used
since the 1950s, starting with the Particle In Cell Method (PIC).
We refer the reader to the classical textbooks~\cite{Birdsall2004,Hockney1988}
for an introduction to the basic concepts and the history of the PIC method.
The review articles~\cite{Dawson1983} and, more recently, \cite{Verboncoeur2005}
discuss advanced aspects of plasma simulations with particle methods.
An obvious strategy for simulation of the particle system is a direct
summation. The force acting on a particle is determined by a summation over
all interaction partners.
Since particles in a plasma interact via long-range Coulomb forces,
an accurate computation of the acceleration of a single particle
requires a summation over all other particles in the plasma.
This results in a quadratic computational complexity,
which is prohibitively expensive with present computer hardware,
even for medium-sized problems.
Therefore, it is key to find approximations to the forces 
which significantly reduce the computational complexity
but, at the same time, preserve their long-range character
and produce consistent results.

Barnes and Hut~\cite{Barnes1986} proposed an approximation
scheme for gravitational problems which they called treecode.
Their idea is to recursively subdivide the particle system
into nested boxes. In three dimensions
each box is split into eight boxes along the Cartesian axes.
The recursive subdivision is embedded into a tree structure,
from now on referred as the cluster tree.
The typical depth of the cluster tree is $\mathcal{O}(\log N_p)$.
The acceleration of a particle $p$ is computed by iterating through the
cluster tree, starting at its root.
The forces between $p$ and all particles in a well-separated 
cluster are replaced by a single force between $p$ 
and a pseudo particle at the centre of mass
of the cluster with mass equal to the total mass of all particles in this cluster.
This generalises easily to electrostatic problems, where the total mass
has to be replaced by the total charge of the cluster.
As the cluster tree has a depth of $\mathcal{O}(\log N_p)$, the numerical
work for the treecode algorithm is $\mathcal{O}(N_p \log N_p)$.
A very similar idea was proposed by Appel \cite{Appel1985} with two major
differences. Firstly, he uses a binary tree, splitting boxes based on
the medians of positions of the particles, and secondly, he avoids
rebuilding the cluster tree after each time step by a merging strategy for
clusters. Again, his algorithm has a complexity of $\mathcal O(N_p \log N_p)$.
Both methods only use the monopole moment of the particle distribution
for the approximation of the forces. This leads to relatively high errors,
especially in the case of non-uniform particle distributions.
However, both methods can be extended to include further terms of the Taylor
expansion. Computations with Taylor expansions up to order $m$
have a complexity of $\mathcal O((m+1)^3 N_p \log N_p)$.
As the error in the far field decays exponentially, $m$ is
chosen as $m \sim |\log \varepsilon| $, where $\varepsilon$ is a predefined
error threshold.

The Fast Multipole Method (FFM), proposed in \cite{Greengard1987}
for two-dimensional problems, and extended in \cite{Greengard1988}
to three-dimensional problems, is also a tree-based method.
In contrast to the treecode discussed above, the FFM uses a Taylor
expansion of the Newton potential in spherical coordinates up to a
given order $m$, a technique well known in electrostatics.
Whereas in the treecode expansions in only one variable were used,
the FFM simultaneously expands the potential in both variables
in the far field. 
Combined with a suitable iteration through the cluster tree,
the numerical cost for the force evaluation is in $\mathcal O ((m+1)^3 N_p)$.
In its first formulation, FFM was restricted to applications with Newton potentials.
Later, it was expanded to general kernels in \cite{Ying2004} and independently,
Of, Steinbach, and Wendland used FFM for the fast solution of boundary integral equations for the Laplace equation \cite{Of2004}
and elastostatics~\cite{Of2005}.

Another important application for the fast evaluation of Coulomb potentials
are molecule dynamics simulations for crystalline structures.
Usually, one neglects boundaries of the crystal and uses periodic boundary conditions for the molecules
and their self-consistent electric field. Those give rise to an infinite sum for the electric potential
which is split into to two rapidly decaying sums. This is known as Ewald summation.
Evaluating the full sum gives a complexity of $\mathcal{O}(N_p^2)$. Darden, York, and Pedersen~\cite{Darden1993}
combined an interpolation scheme with the Fast Fourier Transform to reduce the complexity
to $\mathcal O(N_p \log N_p)$. The approximation error depends on the number of interpolation points.
However, there method is restricted to structured particle distributions, and, more importantly,
to periodic boundary conditions.

In this paper, we present a unified hierarchical framework
for the grid-free simulation of plasma in the electrostatic case in bounded domains with the help
of modern $\mathcal H^2$-matrices.
Both the particle-particle and the particle-boundary interactions
have linear complexity in the number of particles.
We propose the usage of interpolation for the approximation in the far field.
It is very easy to implement, as it only needs the value of a rather
general kernel function at the interpolation points
and furthermore, it is directly applicable to 
the approximation of vector-valued functions.
The contribution of the boundary values to the electric field
are computed via the Boundary Element Method, which only requires a
discretisation of the boundary of the domain. This reduces
the three-dimensional problem posed on the whole domain
to a system of integral equations on a two-dimensional manifold.
Similar ideas have already been presented in \cite{Christlieb2004b, Christlieb2006, Christlieb2008, Christlieb2004a}.
The authors used a treecode-based approximation scheme
with a boundary integral formulation to simulate plasmas in one-
and two-dimensional domains.
Altough theory predicts a complexity of $\mathcal{O}((m+1)^3 N_p \log N_p)$ for their algorithm,
they numerically observe nearly linear scaling in the number of particles.
As we are using $\mathcal H^2$-matrices, we conclude from the theory of hierarchical matrices,
that our algorithm has linear complexity, both in the number of
particles and the number of elements of the surface mesh.
This is supported by our numerical results.
Additionally, we start with the representation formula for the Poisson equation
and systematically approximate the discretised boundary integral operators by $\mathcal H^2$-matrices.
In this way, we treat the particle and the boundary part evenly in terms of the approximation schemes we use.

This article is organised as follows.
\Cref{sec:vlasov-poisson} reviews the Vlasov--Poisson system.
The basic concepts of boundary integral equations and the Boundary
Element Method are given in \cref{sec:BEM}.
In \cref{sec:hierarchical}, we discuss hierarchical approximation
techniques for Nystr\"om and Galerkin matrices.
Important aspects of the implementation of our method in a computer program
are presented in \cref{sec:implementation}.
Numerical examples validating our approach are given in \cref{sec:numerics}.

%% file: vlasov-system.tex
\section{Vlasov--Poisson system}
\label{sec:vlasov-poisson}
If the characteristic velocity $V_0$ of the particle system is small compared to the speed of light $c$,
the dynamics of $N_p$ charged particles with positions $(x_i)_{i=1}^{N_p}$,
velocities $(v_i)_{i=1}^{N_p}$, masses $(m_i)_{i=1}^{N_p}$ and charges $(q_i)_{i=1}^{N_p}$
is given by 
\begin{equation}\label{eq:newton-classic}
\begin{aligned}
	\dot x_i &= v_i, \\
	\dot v_i &= \frac{q_i}{m_i} E(x_i), \quad i=1,\dots,N_p,
\end{aligned}
\end{equation}
coupled with the electrostatic approximation of Maxwell's equations for the electric field $E$
\begin{equation}\label{eq:maxwell-equations-electrostatic}
\begin{aligned}
\operatorname{div} E = \frac{1}{\varepsilon_0} \sum\limits_{j=1}^{N_p} q_j \delta_{x_j},
\quad \operatorname{rot} E = 0,
\end{aligned}
\end{equation}
where $\varepsilon_0$ is the electric constant.

Since $\operatorname{rot} E = 0$, there exists a scalar potential $\phi$ with $E = - \nabla \phi$.
The electrostatic Maxwell's equations can be then expressed as a scalar Poisson equation,
\[
-\Delta \phi = \frac{1}{\varepsilon_0} \sum_{j=1}^{N_p} q_j \delta_{x_j}.
\]
Together with the decay condition for the electric field,
\[
-\nabla \phi(x) = \mathcal{O}\left( \frac{1}{|x|^2} \right), \quad |x| \to \infty,
\]
we obtain the unique solution by applying the Newton potential $N$,
\[
\phi = N \frac{1}{\varepsilon_0} \sum\limits_{j=1}^{N_p} q_j \delta_{x_j}
= \frac{1}{\varepsilon_0} \sum\limits_{j=1}^{N_p} q_j U(\cdot, x_j),
\]
where $U$ is the fundamental solution of the Laplace operator,
\begin{equation}\label{eq:fundamental-solution}
U(x, y) = \frac{1}{4 \pi |x - y|}, \quad x \neq y,
\end{equation}
and the Newton potential for smooth functions with compact support is given by
\begin{equation}\label{eq:newton-potential-free-space}
N \psi (x) = \int_{\mathbb R^3} U(y, x) \psi(y) \, \text dy, \quad x \in \mathbb R^3,
\end{equation}
which is extended to distributions by duality.
The electric field has the form
\begin{equation}\label{eq:electric-field-vaccum}
E = -\nabla \phi = - \frac{1}{\varepsilon_0} \sum\limits_{j=1}^{N_p} q_j \nabla U(\cdot, x_j).
\end{equation}
Plugging~\eqref{eq:electric-field-vaccum} into~\eqref{eq:newton-classic} and excluding self-interactions yields
\begin{equation}\label{eq:newton-electrostatic}
\begin{aligned}
\dot x_i &= v_i, \\
\dot v_i &= -\frac{1}{m_i} \sum\limits_{\substack{j = 1 \\ j \neq i}}^{N_p} \frac{q_i q_j}{\varepsilon_0} \nabla_{x_i} U(x_i, x_j), \quad i=1,\dots,N_p.
\end{aligned}
\end{equation}
This equation is not feasible to describe the time evolution of our plasma as $N_p$ is in the order of $10^{20}$
and therefore out of reach for a direct simulation. We therefore ask for an appropriate limit $N_p \to \infty$
which would give us an easier to handle equation.

One possibility is known as the mean field~\cite{Spohn1991} or pulverisation limit~\cite{Nicholson1983} which is based on a special scaling of charges
and masses.
For simplicity let us assume that our plasma consists of only one species with charge $q_0$ and mass $m_0$.
The masses and charges of the particles are scaled by $1/N_p$,
\[
m_i = m_0 / N_p,~q_i = q_0 / N_p, \quad i=1,\dots,N_p.
\]
This changes~\eqref{eq:newton-electrostatic} to
\begin{equation}\label{eq:newton-mean-field}
\begin{aligned}
	\dot x_i &= v_i, \\
	\dot v_i &= -\frac{1}{N_p} \frac{q_0^2}{m_0 \varepsilon_0} \sum\limits_{\substack{j = 1 \\ j \neq i}}^{N_p} \nabla_{x_i} U(x_i, x_j), \quad i=1,\dots,N_p.
\end{aligned}
\end{equation}
In their pioneering work, Neunzert and Wick~\cite{Neunzert1972,Neunzert1973,Neunzert1974} (for an English version of their ideas see~\cite{Neunzert1984} and
also~\cite{Spohn1991} for a historic review)
show the convergence of~\eqref{eq:newton-mean-field} for $N_p \to \infty$ to the solution $f: (0, \infty) \times \mathbb R^3 \times \mathbb R^3 \to [0, \infty)$
of the Vlasov-Poisson system,
\begin{equation}
\begin{gathered}\label{eq:vlasov-poisson-free-space}
\partial_t f + v \cdot \nabla_x f + \frac{q_0}{m_0} E \cdot \nabla_v f = 0, \\
E = - \nabla \phi, \\
-\Delta \phi = \frac{q_0}{\varepsilon_0} \int_{\mathbb R^3} f \, \text dv, \quad -\nabla \phi(x) = \mathcal{O}\left( \frac{1}{|x|^2} \right), \quad |x| \to \infty,
\end{gathered}
\end{equation}
in the weak-$\ast$ topology of measures.
They rely on a regularisation $G_\varepsilon$ of $G = -\nabla_x U$ which is assumed to be continuous and bounded.
Possible choices for $G_\varepsilon$ include a mollified version of $G$ or the gradient of
\[
U_\varepsilon(x,y) = \frac{1}{4\pi |x-y| + \varepsilon}, \quad x,y \in \mathbb R^3.
\]
The parameter $\varepsilon$ is time-dependent and tends to $0$ as $t \to \infty$, see~\cite{Ganguly1989, Wollman2000} and the references cited therein.
In a recent work, Lazarovici and Pickl~\cite{Lazarovici2017} show convergence to the Vlasov--Poisson system for a scaling that only depends on $N_p$, $\varepsilon = N_p^{-1/3 + o(1)}$.
This is nearly optimal in the sense that the mean distance between two particles scales like $N_p^{-1/3}$.
Their regularisation of $G$ has to satisfy three conditions,
\begin{enumerate}
	\item $\exists c_1 > 0 ~ \forall x,y \in \mathbb R^3, x \neq y: |G_\varepsilon(x,y)| \leq c_1 / |x-y|^2,~ |\nabla_x G_\varepsilon(x,y)| \leq c_1 / |x-y|^3$,
	\item $\forall x,y \in \mathbb R^3, |x-y| \geq \varepsilon: G_\varepsilon(x,y) = G(x,y)$,
	\item $\exists c_2 > 0 ~ \forall x,y \in \mathbb R^3, |x-y| < \varepsilon: |G_\varepsilon(x,y)| \leq c_2/\varepsilon^2, ~ |\nabla_x G_\varepsilon(x,y)| \leq c_2/\varepsilon^3$.
\end{enumerate}
In order to systematically derive suitable regularisations for the interaction force $G$, we can also regularise
the charge distribution. Applying the Newton potential~\eqref{eq:newton-potential-free-space} then gives a regularisation
which is consistent with the charge density in the Poisson equation. 
Furthermore, this allows us to apply the standard theory of Sobolev spaces for elliptic problems.
For our implementation we choose a radial step function,
\[
\delta^{\varepsilon}_{y} = \frac{1}{|B_{\varepsilon}(y)|} \mathbbm{1}_{B_{\varepsilon}(y)}, \quad y \in \mathbb R^3,
\]
for which we have
\begin{equation}\label{eq:potential-regularisation}
U_\varepsilon(x, y) = N \delta^{\varepsilon}_{y}(x) = \frac{1}{4 \pi} 
\begin{cases}
	\dfrac{3}{2 \varepsilon} - \dfrac{|x-y|^2}{2 \varepsilon^3}, & |x-y| < \varepsilon \\[0.5em]
	\dfrac{1}{|x - y|}, & |x-y| \geq \varepsilon
\end{cases},
\quad x, y \in \mathbb R^3.
\end{equation}
Applying the gradient to~\eqref{eq:potential-regularisation} yields
\begin{equation}\label{eq:potential-gradient-regularisation}
G_\varepsilon(x,y) = \frac{1}{4 \pi} 
\begin{cases}
	\dfrac{1}{\varepsilon^3}(x-y), & |x-y| < \varepsilon \\[0.5em]
	\dfrac{x-y}{|x - y|^3}, & |x-y| \geq \varepsilon
\end{cases},
\quad x, y \in \mathbb R^3.
\end{equation}
It is easy to check that $G_\varepsilon$ is a bounded Lipschitz continuous function that satisfies the aforementioned conditions
on the regularisation. Furthermore, $\delta^\varepsilon_y \in L_2(\mathbb R^3)$ which simplifies the analysis in \cref{sec:BEM}
when working with trace operators.

In this paper, we consider the Vlasov-Poisson system~\eqref{eq:vlasov-poisson-free-space} in a bounded domain~$\Omega$.
Instead of decay conditions on the potential $\phi$ we now prescribe Dirichlet or Neumann conditions on the boundary $\partial \Omega$.
Additionally, we also need boundary conditions for the distribution function $f$. We primarily choose absorption,
i.e. $f=0$ on $\partial \Omega$. In its nondimensional form, the Vlasov-Poisson system reads
\begin{equation}\label{eq:vlasov-poisson}
\begin{gathered}
\partial_t f(t,x,v) + v \cdot \nabla_x f(t,x,v) + E(t,x) \cdot \nabla_v f(t,x,v) = 0, \\
E(t,x) = - \nabla_x \phi(t,x), \\
-\Delta_x \phi(t,x) = \frac{1}{\beta} \int_{\mathbb R^3} f(t,x,v) \, \text dv,
\end{gathered}
\end{equation}
for $(t,x,v) \in (0, \infty) \times \Omega \times \mathbb R^3$.
Here, $\beta = \left(\lambda_D / L_0\right)^2$ is the square of the non-dimensional quotient of the Debye length
\[
\lambda_D = \sqrt{\frac{\varepsilon_0 k_B T_0}{n_0 q_0^2}}
\]
and the characteristic length $L_0$ of $\Omega$.
Furthermore, $k_B$ is the Boltzmann constant and $T_0, n_0, q_0$
denote the characteristic temperature, particle density and
charge of the plasma, respectively.
Initial conditions for $f$ are usually linear combinations of Maxwellians,
\[
M_{\rho_m, V, T}(x,v) = \frac{\rho_m(x)}{(2 \pi T(x))^{3/2}} \exp \left( -\frac{|v - V(x)|^2}{2 T(x)} \right), \quad (x, v) \in \Omega \times \mathbb R^3,
\]
where $\rho_m$ is the mass density, $V$ is the macroscopic bulk velocity, and $T$ the temperature distribution inside the plasma.

For the numerical treatment of~\eqref{eq:vlasov-poisson}, we sample $f$ by $N_p$ macroparticles,
\begin{equation}\label{eq:distribution-function-measures}
f(t,\cdot, \cdot) \approx \frac{|\Omega|}{N_p} \sum\limits_{j=1}^{N_p} \delta_{x_j(t)} \, \delta_{v_j(t)}, \quad t > 0
\end{equation}
and regularise the charge density
\begin{equation}\label{eq:regularised-charge-density}
\rho(t,x) = q_0 \int_{\mathbb R^3} f(t,x,v) \, \text dv
\approx q_0 \frac{|\Omega|}{N_p} \sum\limits_{j=1}^{N_p} \frac{1}{|B_\varepsilon(x_j(t))|} \mathbbm{1}_{B_\varepsilon(x_j(t))},
\end{equation}
where $|\Omega|$ is the volume of the (rescaled) domain.
The Vlasov equation for the approximation~\eqref{eq:distribution-function-measures} is equivalent to a system of ODEs,
\[
\begin{aligned}
\dot x_i &= v_i, \\
\dot v_i &= - \frac{q_0}{m_0} \nabla \phi(x_i) \quad i=1,\dots,N_p,
\end{aligned}
\]
where $\phi$ is the solution to the boundary value problem
\begin{equation}\label{eq:dirichlet-problem}
\begin{aligned}
-\Delta \phi &= \frac{q_0}{\beta} \frac{|\Omega|}{N_p} \sum\limits_{j=1}^{N_p} \frac{1}{|B_\varepsilon(x_j(t))|} \mathbbm{1}_{B_\varepsilon(x_j(t))} & & \mbox{in  $\Omega$,} \\
\phi &= g_D & & \mbox{on $\Gamma = \partial \Omega$},
\end{aligned}
\end{equation}
and we assume a Dirichlet problem for simplicity.
Keeping in mind that a special solution to above equation is given by the Newton potential~\eqref{eq:newton-potential-free-space}
a particular solution $\phi_p$ of the Poisson equation above for a
fixed time $t > 0$ is
\begin{equation}\label{eq:phi_p}
\phi_p(t,x) = 
\frac{q_0}{\beta} \frac{|\Omega|}{N_p} \sum_{j=1}^{N_p} U_\varepsilon(x, x_j(t)), \quad x \in \Omega.
\end{equation}

In order to find a solution of the BVP~\eqref{eq:dirichlet-problem}
with the help of $\phi_p$, we have
to solve the auxiliary problem
\begin{equation}\label{eq:aux-dirichlet-problem}
\begin{aligned}
-\Delta \phi_0 &= 0  & & \text{in } \Omega, \\
\phi_0 &= g_D-\phi_p & & \text{on } \Gamma = \partial \Omega.
\end{aligned}
\end{equation}
The solution of the original problem is now
\begin{equation}\label{eq:decomposition-phi}
\phi = \phi_0 + \phi_p,
\end{equation}
and the electric field at the time $t$ in the position of particle $i$
is computed as
\begin{equation}\label{eq:electric-field}
E(t,x_i(t)) = -\nabla \phi_0(t,x_i(t)) + \frac{q_0}{\beta}
\frac{|\Omega|}{N_p}\sum\limits_{\substack{j=1\\ j \neq i}}^{N_p} 
	G_\varepsilon(x_i(t), x_j(t)).
\end{equation}
Note that this representation is consistent with the mean field scaling~\eqref{eq:newton-mean-field}
that leads to the Vlasov-Poisson equation. For homogeneous boundary conditions, $\phi_0$ scales
like $1/N_p$. This is also true for a pure Neumann or a mixed boundary value problem.

Whereas the evaluation of $\phi_p$ is grid-free by its nature,
the numerical treatment of equation~\eqref{eq:aux-dirichlet-problem}
involves, as a rule, the discretisation of the domain.
For simple domains of toroidal or rectangular shape the discretisation of the Poisson equation
on structured grids leads to linear systems whose solutions are usually found by means of the Fast Fourier Transform~\cite{Birdsall2004}.
In order to evaluate the electric field at the positions of the particles and to couple the charge density
with the grid, a regularisation of $\rho_\text{total}$ is needed. The electric field at the positions of the particles
is obtained by interpolation from the grid nodes.
Structured meshes also work for complex domains, where one relies on so-called cut cells near the boundary,
see~\cite{Kolobov2016} for the electrostatic case and~\cite{Nieter2009} for the full Maxwell system.
Without a suitable post-processing, degenerated cut cells with small side lengths put an additional constrain on the CFL
condition for an explicit scheme~\cite{Kolobov2016,Nieter2009}.
Contrarily, for the proposed Boundary Element Methods we use to solve~\eqref{eq:aux-dirichlet-problem},
no volume discretisation is needed and therefore we introduce no further restriction on the CFL condition.
By the use of the representation formula, we can compute the electric field at each given point inside the domain,
i.\,e. the positions of the particles.
Note that this decouples the particle discretisation of the distribution function and the discretisation of the
Poisson equation for the electric potential.
The number of particles does not affect the accuracy of the electric field which is only controlled by the mesh size of the boundary mesh.
The particles move freely through the volume $\Omega$.
As written in~\eqref{eq:electric-field}, the electric field is split into two parts.
The free space interaction of the particles which ignores boundary conditions and a correction term which solves~\eqref{eq:aux-dirichlet-problem}
and depends on boundary conditions.
For this, we propose Boundary Element Methods.
The accuracy of the electric field only depends on the error made in the approximation of $-\nabla \phi_0$ but not on the number of particles.
In contrast to PIC methods, there is no rule of thumb connecting the number of particles and the number of triangles of the surface mesh.
In the following sections, we first review Boundary Elements Methods and the discretisation for mixed problems.
With the notation we define there, we are able to give a first formulation of our algorithm with quadratic complexity in \cref{sec:BEM-Vlasov-Poisson}.
Afterwards, we discuss how to accelerate it and to reduce the complexity from quadratic to linear in the number of particles and the number of triangles.

%% file: bem.tex
\section{Boundary Element Method}
\label{sec:BEM}
The Boundary Element Method (BEM) is reviewed for the general Poisson problem with mixed boundary conditions on a bounded polyhedral domain $\Omega\subset\mathbb R^3$ with boundary $\Gamma=\partial\Omega$. Furthermore, $\Gamma=\overline\Gamma_D\cup\overline\Gamma_N$ is split into an open Dirichlet part $\Gamma_D$ and an open Neumann part $\Gamma_N$.

We first define the abstract mathematical framework for BEM based on fractional Sobolev spaces on the boundary. The section concludes with the definition of potential operators.
The reader interessed in the discretisation of boundary integral equations may skip this section and start with the section on Galerkin discretisation.
\subsection{Boundary integral equations}
The classical theory for boundary integral equations is based on square-integrable functions and their (weak) derivatives.
Let $L_2(\Omega)$ denote the space of square-integrable functions,
\[
L_2(\Omega) = \left\{f: \Omega \to \mathbb R: \int_{\Omega} |f(x)|^2 \, \text dx < \infty \right\},
\]
which is a Hilbert space with respect to the inner product
\[
(f,g)_{L_2(\Omega)} = \int_\Omega f(x) g(x) \, \text dx, \quad f,g \in L_2(\Omega).
\]
Analogously, $L_2(\Gamma)$, the space of square-integrable functions on the boundary, is defined.
The Sobolev space $H^1(\Omega)$ consists of functions in $L_2(\Omega)$ which have a weak gradient in $L_2(\Omega)$, i.e.
for $f \in H^1(\Omega)$ there exists a $g \in L_2(\Omega)^3$ such that for all $\varphi \in C_0^\infty(\Omega)^3$
\[
\int_\Omega g(x) \cdot \varphi(x) \, \text dx = - \int_\Omega f(x) \nabla \cdot \varphi(x) \, \text dx,
\]
where $C_0^\infty(\Omega)$ denotes the set of infinitely often differentiable functions with compact support in $\Omega$.
Equipped with the inner product
\[
(f, g)_{H^1(\Omega)} = (f, g)_{L_2(\Omega)} + (\nabla f, \nabla g)_{L_2(\Omega)^3},
\]
the Sobolev space $H^1(\Omega)$ is a Hilbert space.
Similar spaces, called Sobolev--Slobodekii or fractional Sobolev spaces, can be defined on the boundary~\cite{Steinbach2007},
\[
H^{1/2}(\Gamma) = \left\{ f \in L_2(\Gamma): \int_\Gamma \int_\Gamma \frac{|f(x)-f(y)|^2}{|x-y|^{3}} \, \text ds_y \, \text ds_x < \infty \right\},
\]
which forms a Hilbert spaces with inner product
\[
(f, g)_{H^{1/2}(\Gamma)} = (f, g)_{L_2(\Gamma)} + \int_\Gamma \int_\Gamma \frac{(f(x)-f(y))(g(x)-g(y))}{|x-y|^{3}} \, \text ds_y \, \text ds_x,~ f,g \in H^{1/2}(\Gamma).
\]
Spaces with negative indices like $H^{-1/2}(\Gamma)$ are defined as the dual spaces with respect to $\langle \cdot, \cdot \rangle_{\Gamma}$, the extension of the $L_2(\Gamma)$-inner product.
The norm in $H^{-1/2}(\Gamma)$ is given by
\[
\| u \|_{H^{-1/2}(\Gamma)} = \sup\limits_{f \in H^{1/2}(\Gamma)\setminus \{0\}} \frac{\langle f, u \rangle_\Gamma}{\| f \|_{H^{1/2}(\Gamma)}}, \quad u \in H^{-1/2}(\Gamma).
\]
For example, $H^{1/2}(\Gamma)$ contains continuous functions but not piecewise continuous functions with discontinuous jumps. But those functions are included
in $H^{-1/2}(\Gamma)$. This is an important observation for choosing ansatz and test functions for the Galerkin formulation later.
For a mixed formulation, we also need fractional Sobolev spaces on open subsets $\Gamma_1$ of the boundary.
Then, $H^{1/2}(\Gamma_1)$ is defined by restrictions of functions from $H^{1/2}(\Gamma)$,
\[
H^{1/2}(\Gamma_1) = \{ f\left|_{\Gamma_1}\right.: f \in H^{1/2}(\Gamma_1) \},
\]
with norm
\[
\| g \|_{H^{1/2}(\Gamma_1)} = \inf\limits_{f\left|_{\Gamma_1}\right. = g} \| f \|_{H^{1/2}(\Gamma_1)}, \quad g \in H^{1/2}(\Gamma_1).
\]
The space $H^{-1/2}(\Gamma_1)$ is formed by all continuous linear functionals acting on functions in $H^{1/2}(\Gamma)$ with support in $\Gamma_1$.
Note that duality is understood with respect to $\langle \cdot, \cdot, \rangle_\Gamma$.

Given a volume source term $g_V\in L_2(\Omega)$, a Dirichlet datum $g_D\in H^{1/2}(\Gamma_D)$ as well as a Neumann datum $g_N\in H^{-1/2}(\Gamma_N)$, the Poisson problem reads
\begin{equation}\label{eq:BEM-PoissonProblem}
 \begin{aligned}
  -\Delta \phi & = g_V & &\mbox{in } \Omega,\\
  \phi & = g_D & & \mbox{on } \Gamma_D,\\
  n_\Omega\cdot\nabla \phi & = g_N & & \mbox{on } \Gamma_N,
 \end{aligned}
\end{equation}
where $n_\Omega$ denotes the outward unit normal vector on~$\Gamma$. The boundary value problem is considered in the weak sense, such that the solution is sought in the Sobolev space~$H^1(\Omega)$. We may follow the idea of the previous section and construct a particular solution~$\phi_p$ in order to homogenise the right hand side of the differential equation. An appropriate choice is the Newton potential
\begin{equation}\label{eq:BEM-phi_p}
 \phi_p(x) = (Ng_V)(x) = \int_\Omega U(x,y)\,g_V(y)\,\text{d}y
 \quad\mbox{for } x\in\mathbb R^3,
\end{equation}
where $U(x,y)$ is the fundamental solution given in~\eqref{eq:fundamental-solution}. For $g_V=\rho_\text{total}/\beta$ we recover~\eqref{eq:phi_p}. The problem~\eqref{eq:BEM-PoissonProblem} has a unique solution that admits for $x\in\Omega$ the representation formula
\begin{equation}\label{eq:BEM-representation-formula}
 \phi(x)=\int_{\Gamma}U(x,y)\gamma_1\phi(y)\,\text{d}s_y - \int_{\Gamma}\gamma_{1,y} U(x,y) \gamma_0\phi(y) \,\text{d}s_y + (Ng_V)(x)
\end{equation}
where $\gamma_0\phi$ denotes the Dirichlet and $\gamma_1\phi$ the Neumann trace of the unknown solution~$\phi$. For sufficiently smooth data and $x\in\Gamma$ it holds
\[
 \gamma_0\phi(x) = \phi\big|_{\Gamma}(x)
 \qquad\mbox{and}\qquad
 \gamma_1\phi(x) = \lim_{\Omega\ni \tilde x\to x} n_\Omega\cdot\nabla\phi(\tilde x).
\]
These trace operators can be extended to linear bounded operators with the following mapping properties~\cite{McLean2000}:
\[
 \gamma_0:H^1(\Omega)\to H^{1/2}(\Gamma)
 \qquad\mbox{and}\qquad
 \gamma_1: H_{\Delta}^1(\Omega) \to H^{-1/2}(\Gamma),
\]
where $f \in H_{\Delta}^1(\Omega)$ iff $f \in H^1(\Omega)$ and $\Delta f \in L_2(\Omega)$.
We apply the trace operators to the representation formula~\eqref{eq:BEM-representation-formula} and obtain the system of equations
\begin{equation}\label{eq:BEM-System}
 \begin{pmatrix} \gamma_0\phi \\ \gamma_1\phi \end{pmatrix} =
 \begin{pmatrix} \frac{1}{2} I- K & V \\ W & \frac{1}{2}I+K^\prime \end{pmatrix}  \begin{pmatrix} \gamma_0\phi \\ \gamma_1\phi \end{pmatrix} +
 \begin{pmatrix} N_0 g_V \\ N_1 g_V \end{pmatrix}.
\end{equation}
This system contains the standard boundary integral operators which are well studied, see, e.g.,~\cite{McLean2000,SauterSchwab2011,Steinbach2007}. For $x\in\Gamma$, we have the single-layer potential operator
\[
 (V \zeta)(x) = \gamma_0\int_{\Gamma}\!U(x,y)\zeta(y)\,\text{d}s_{y} \quad\mbox{for }\zeta\in H^{-1/2}(\Gamma),
\]
the double-layer potential operator
\[
 (K \xi)(x) = \lim_{\varepsilon\to 0}\int\limits_{y\in\Gamma:\|y-x\|\geq\varepsilon}\hspace{-0.4cm}\gamma_{1,y} U(x,y)\xi(y)\,\text{d}s_{y} \quad\mbox{for } \xi\in H^{1/2}(\Gamma),
\]
where $\gamma_{1,y}$ means that the Neumann trace $\gamma_1$ only acts on the $y$-variable, and the adjoint double-layer potential operator
\[
 (K^\prime \zeta)(x) = \lim_{\varepsilon\to 0}\int\limits_{y\in\Gamma:\|y-x\|\geq\varepsilon}\hspace{-0.4cm}\gamma_{1,x} U(x,y)\zeta(y)\,\text{d}s_{y} \quad\mbox{for } \zeta\in H^{-1/2}(\Gamma),
\]
as well as the hypersingular integral operator
\[
 (W \xi)(x) = -\gamma_1\int_{\Gamma}\gamma_{1,y} U(x,y)\xi(y)\,\text{d}s_{y} \quad\mbox{for } \xi\in H^{1/2}(\Gamma),
\]
and $N_0g_V = \gamma_0 Ng_V$ as well as $N_1g_V = \gamma_1 Ng_V$.

\subsection{Galkerin discretisation}\label{sec:BEM-Galkerin}
Obviously, if the traces $\gamma_0\phi$ and $\gamma_1\phi$ of the unknown solution~$\phi$ are known, the representation formula~\eqref{eq:BEM-representation-formula} can be used to evaluate~$\phi$ inside the domain~$\Omega$. However, these traces are only known on parts of the boundary according to~\eqref{eq:BEM-PoissonProblem}. Thus, we aim to approximate them on the whole boundary~$\Gamma$ with the help of a Galerkin BEM, following~\cite{Steinbach2007}. Therefore, let $\Gamma$ be meshed by a quasi-uniform, conforming surface triangulation $\mathcal T = \left\{ \Gamma_k \right\}_{k=1}^{N_\Gamma}$ that is shape-regular in the sense of Ciarlet with $N_\Gamma$ triangles and $M_\Gamma$ nodes. We apply the conforming approximation spaces
\[
 S_h^0(\Gamma) = \Span\big\{\varphi_k^0\big\}_{k=1}^{N_\Gamma}\subset H^{-1/2}(\Gamma),
 \quad\mbox{and}\quad
 S_h^1(\Gamma) = \Span\big\{\varphi_i^1\big\}_{i=1}^{M_\Gamma}\subset H^{1/2}(\Gamma),
\]
where $\varphi_k^0$ denotes the piecewise constant function that is one on the triangle of index~$k$ and zero else, and $\varphi_i^1$ denotes the usual hat function corresponding to the node with index~$i$. For simplicity, we write $\phi=\gamma_0\phi$ and assume that the triangles and nodes are numbered in such a way that the triangles for $k=1,\ldots,N_D$ lie in~$\Gamma_D$ and the nodes for $i=1,\ldots,M_N$ are the ones without Dirichlet condition. We seek the approximation of the Dirichlet trace as 
\begin{equation}\label{eq:BEM-Ansatz-Dirichlet}
 \phi_h(x) 
 = \phi_{N,h}(x) + \phi_{D,h}(x) 
 = \sum_{i=1}^{M_N}\phi_i\varphi_i^1(x) + \sum_{i=M_N+1}^{M_\Gamma}\phi_i\varphi_i^1(x)
\end{equation}
and the Neumann trace as
\begin{equation}\label{eq:BEM-Ansatz-Neumann}
 t_h(x) 
 = t_{D,h}(x) + t_{N,h}(x) 
 = \sum_{k=1}^{N_D}t_k\varphi_k^0(x) + \sum_{k=N_D+1}^{N_\Gamma}t_k\varphi_k^0(x)
\end{equation}
with vectors $\underline\phi_{N,h}=(\phi_1,\ldots,\phi_{M_N})^\top\in\mathbb R^{M_N}$ and $\underline t_{D,h}=(t_1,\ldots,t_{N_D})^\top\in\mathbb R^{N_D}$, respectively, and $\underline\phi_{D,h}$ and $\underline t_{N,h}$ accordingly. 
The coefficients $\phi_i$, $i=M_N+1,\ldots,M_\Gamma$ and $t_k$, $k=N_D+1,\ldots,N_\Gamma$ are determined by interpolation of the given boundary data in~\eqref{eq:BEM-PoissonProblem}. 
Inserting the ansatz~\eqref{eq:BEM-Ansatz-Dirichlet} and~\eqref{eq:BEM-Ansatz-Neumann} into~\eqref{eq:BEM-System}, testing with 
$\varphi^0_k$, $k=1,\ldots,N_D$ 
and 
$\varphi^1_i$, $i=1,\ldots,M_N$, 
respectively, and integrating over~$\Gamma$ yields
\begin{multline}\label{eq:BEM-block-system}
 \begin{pmatrix}
  V_h^{DD} & -K_h^{DN} \\ {K_h^{DN}}^\top & W_h^{NN}
 \end{pmatrix}
 \begin{pmatrix}
  \underline t_{D,h} \\ \underline\phi_{N,h}
 \end{pmatrix} \\
 = 
 \begin{pmatrix}
  \frac{1}{2} M_h^{DD}+ K_h^{DD} & -V_h^{DN} \\ -W_h^{ND} & \frac{1}{2}{M_h^{NN}}^\top-{K_h^{NN}}^\top 
 \end{pmatrix}
 \begin{pmatrix}
  \underline\phi_{D,h} \\ \underline t_{N,h} 
 \end{pmatrix}
 -
 \begin{pmatrix}
  \underline N_0^{D} \\ \underline N_1^{N} 
 \end{pmatrix}.
\end{multline}
The matrices are defined by
\begin{equation}
 \begin{aligned}\label{eq:galerkin-bem-matrices}
 V_h[\ell,k] &= (V\varphi^0_k,\varphi^0_\ell)_{L_2(\Gamma)}, &
 W_h[j,i] &= (D\varphi^1_i,\varphi^1_j)_{L_2(\Gamma)}, \\
 K_h[\ell,i] &= (K\varphi^1_i,\varphi^0_\ell)_{L_2(\Gamma)}, &
 M_h[\ell,i] &= (\varphi^1_i,\varphi^0_\ell)_{L_2(\Gamma)},
 \end{aligned}
\end{equation}
where $i,j=1,\ldots,M_\Gamma$ and $k,\ell=1,\ldots,N_\Gamma$,
with the block structure 
\begin{equation}\label{eq:discrete-potentials}
 \begin{aligned}
 V_h &= \begin{pmatrix} V_h^{DD} & V_h^{DN} \\ V_h^{ND} & V_h^{NN}\end{pmatrix}, &
 W_h &= \begin{pmatrix} W_h^{NN} & W_h^{ND} \\ W_h^{DN} & W_h^{DD}\end{pmatrix}, \\
 K_h &= \begin{pmatrix} K_h^{DN} & K_h^{DD} \\ K_h^{NN} & K_h^{ND}\end{pmatrix}, &
 M_h &= \begin{pmatrix} M_h^{DN} & M_h^{DD} \\ M_h^{NN} & M_h^{ND}\end{pmatrix},
 \end{aligned}
\end{equation}
representing the Dirichlet and Neumann boundary parts of the matrices.
Fully written out, the entries for the single- and double-layer potential read
\[
V_h[\ell, k] = \frac{1}{4 \pi} \int_{\Gamma_\ell} \int_{\Gamma_k} \frac{1}{|x - y|} \, \text ds_y \, \text ds_x,
\]
and
\[
K_h[\ell, i] = \frac{1}{4 \pi} \int\limits_{\Gamma_\ell} \int\limits_{\operatorname{supp} \varphi^1_i} \frac{(x-y) \cdot n_\Omega(y)}{|x-y|^3} \varphi^1_i(y) \, \text ds_y \, \text ds_x,
\]
for $k,\ell=1,\dots,N_\Gamma$ and $i=1,\dots,M_\Gamma$. Note that these integrals are singular if the supports of ansatz and test functions overlap.
Therefore, special quadratures are used to accurately compute the matrix entries. The most general method, which is sometimes called black box quadrature,
was developed by Sauter and Schwab, see~\cite{SauterSchwab2011}. It is based on a suitable regularisation of the kernel function and utilises tensorised 
Gau\ss-Legendre quadrature on $[-1,1]^{4}$. The quadrature error is proven to decay exponentially with the number of quadrature points.

Furthermore, we used
\[
  \underline N_0[\ell] = (N_0g_V,\varphi^0_\ell)_{L_2(\Gamma)}
  \quad\mbox{with}\quad
  \underline N_0 = (\underline N^{D}_0, \underline N^{N}_0)^\top
\]
for $\ell=1,\ldots,N_\Gamma$ and
\[
  \underline N_1[j] = (N_1g_V,\varphi^1_j)_{L_2(\Gamma)}
  \quad\mbox{with }
  \underline N_1 = (\underline N^{N}_1, \underline N^{D}_1)^\top
\]
for $j=1,\ldots,M_\Gamma$.
Since in our case $N_0g_V$ is computed easily using \eqref{eq:BEM-phi_p}, we exploit the identity
\begin{equation}\label{eq:neumann-newton-potential}
  N_1g_V = \left(-\tfrac12I+K^\prime\right)V^{-1}N_0g_V
\end{equation}
in order to approximate $\underline N_1$ and to avoid volume integrals. We refer the interested reader to~\cite{OfSteinbachUrthaler2010} for more details.

For a pure Dirichlet problem, i.e.\ $\Gamma_N=\varnothing$, the system reduces to
\begin{equation}\label{eq:BIE-pure-Dirichlet}
  V_h\underline t_h = \left(\tfrac12M_h+K_h\right)\underline\phi_h - \underline N_0.
\end{equation}
We can omit the Newton potential when utilising the proposed decomposition~\eqref{eq:decomposition-phi} with $\phi_0$ as solution of~\eqref{eq:aux-dirichlet-problem}. This ansatz yields for the approximation of the Neumann trace $t_{0,h}\approx \gamma_1\phi_0$ the system of linear equations
\[
  V_h\underline t_{0,h} = \left(\tfrac12M_h+K_h\right)\underline\phi_{0,h},
\]
where $\phi_{0,h}\approx\gamma_0\phi_0=g_D-\gamma_0u_p$.

For a pure Neumann problem, i.e.\ $\Gamma_D=\varnothing$, the system also reduces. The hypersingular integral operator, however, is not invertible on $H^{1/2}(\Gamma)$ and thus, the stabilised system~\cite{Steinbach2007}
\begin{equation}\label{eq:BIE_LGSNeumannDirichletVF}
 \widetilde{W}_h\underline \phi_h = \left(\tfrac{1}{2}M_h^\top-K_h^\top\right)\underline t_h - \underline N_1,
\end{equation}
is considered, where
\[
 \widetilde{W}_h = W_h+\alpha\, d_h\,d_h^\top
 \quad\mbox{with}\quad
 d_h[i]=\left(\varphi_i^1,1\right)_{L^2(\Gamma)},
\]
and stabilisation parameter $\alpha>0$.
The system~\eqref{eq:BIE_LGSNeumannDirichletVF} is uniquely solvable since the matrix $\widetilde{W}_h$ is symmetric and positive definite due to the properties of the integral operator $W$. Furthermore the stabilisation ensures that 
\[
 \int_\Gamma t_h(x)\,\text{d}s_x = 0.
\]
At the end of this section, we shortly discuss the approximation error of the Galerkin method.
Since we are only interested in point values of the solution in the interior of the domain, we focus on pointwise error estimates.
We cite the main results and do not give all necessary conditions for the following theorems to hold.
The reader is refered to~\cite{SauterSchwab2011,Steinbach2007} for more details.
\begin{lemma}\label{lem:BEM-error-dirichlet-neumann}
For the numerical approximations $\phi_h$~\eqref{eq:BEM-Ansatz-Dirichlet} and $t_h$~\eqref{eq:BEM-Ansatz-Neumann} the error estimates
\[
\| t_h - \gamma_1 \phi \|_{L_2(\Gamma)} = \mathcal O(h), \quad \| \phi_h - \gamma_0 \phi \|_{L_2(\Gamma)} = \mathcal O(h^2)
\]
hold, where $h$ is the mesh size of $\mathcal T$.
\end{lemma}
From \cref{lem:BEM-error-dirichlet-neumann} pointwise error estimates follow.
\begin{lemma}\label{lem:BEM-pointwise-error}
For $x \in \Omega$ there are $C_1, C_2 > 0$ such that
\[
|\phi(x) - \widetilde \phi_h(x) | \leq C_1 h^3 \quad |\nabla \phi(x) - \nabla \widetilde \phi_h(x) | \leq C_2 h^3
\]
for the pure Dirichlet or Neumann problem. For the mixed problem we have at least quadratic convergence.
Here, $\widetilde{\phi}_h$ denotes the function obtained by plugging $\phi_h$ and $t_h$ into the representation formula~\eqref{eq:BEM-representation-formula}.
\end{lemma}
\subsection{Application to the Vlasov-Poisson system}\label{sec:BEM-Vlasov-Poisson}
We now specify the choice of the data for the general Poisson equation~\eqref{eq:BEM-PoissonProblem} for the self-consistent field of the particles
and give a first naive version of our algorithm for the computation of the electric field in \cref{alg:electric-field-quadratic}.
The volume term $g_V$ is proportional to $\rho$ from~\eqref{eq:regularised-charge-density}, i.e.
\[
g_V = \frac{q_0}{\beta} \frac{|\Omega|}{N_p} \sum\limits_{j=1}^{N_p} \frac{1}{|B_\varepsilon(x_j|} \mathbbm{1}_{B_\varepsilon(x_j}
\]
and therefore
\[
N_0 g_V = \frac{q_0}{\beta} \frac{|\Omega|}{N_p} \sum\limits_{i=1}^{N_p}  U_\varepsilon(\cdot, x_i)
\]
with $U_\varepsilon$ given by~\eqref{eq:potential-regularisation}.
The discretised Dirichlet trace of the Newton potential now is
\[
\underline{N}_0[\ell] = \frac{q_0}{\beta} \frac{|\Omega|}{N_p} 
\sum\limits_{j=1}^{N_p} \left(\int_{\Gamma_\ell} U_\varepsilon(y, x_j) \, \text ds_y \right) \quad \ell=1,\dots,N_\Gamma.
\]
For an efficient implementation into a computer program and for a hierarchical approximation of the electric field,
it is important to restate all operations as matrix-vector products.
We begin with $\underline N_0$ whose computation is expressed as
\begin{equation}\label{eq:newton-dirichlet-mvm}
\underline N_0 = \frac{1}{\beta} \Phi_{\mathcal T} \, w_q,
\end{equation}
where $w_q = (q_0 |\Omega|/N_p)_{i=1}^{N_p}$ is the vector of weighted charges and the entries of $\Phi_{\mathcal T}$ are given by
\[
\Phi_{\mathcal T}[\ell, j] = \int_{\Gamma_\ell} U_\varepsilon(y, x_j) \, \text ds_y,
\]
for $\ell=1,\dots,N_\Gamma$ and $i=1,\dots,N_p$.
In a similar way we reformulate the gradient of the representation formula~\eqref{eq:BEM-representation-formula}.
For this, we define $F_k \in \mathbb R^{N_p \times N_p}$, $S_k \in \mathbb R^{N_p \times N_\Gamma}$, $D_k \in \mathbb R^{N_p \times M_\Gamma}$ for $k=1,2,3$:
\begin{equation}\label{eq:field-boundary-integrals}
S_k[i, \ell] = \int_{\Gamma_\ell} \frac{\partial}{\partial x_i^{(k)}} U(x_i, y) \, \text ds_y, \quad
D_k[i, j] = \int\limits_{\operatorname{supp} \varphi_j} \frac{\partial}{\partial x_i^{(k)}} \gamma_{1,y} U(x_i, y) \varphi_j(y) \, \text ds_y,
\end{equation}
for $i=1,\dots,N_p$, $\ell=1,\dots,N_\Gamma$ and $j=1,\dots,M_\Gamma$.
Furthermore,
\begin{equation}\label{eq:field-free-space}
F_k[i,j] = \frac{\partial}{\partial x_i^{(k)}} U_\varepsilon(x_i, x_j)
\end{equation}
for $i,j=1,\dots,N_p$.
The matrices $S_k$ and $D_k$, $k=1,2,3$ are the contributions of the single- and double-layer potentials to the gradient of the solution, respectively.
The matrices $(F_k)_{k=1}^3$ represent the gradient of the Newton potential, i.e. the free space field of the particles.
We collect the values of the electric field evaluated at the positions of the particles in three vectors,
\[
\underline E_k[i] = E_k(x_i), \quad k=1,2,3,~i=1,\dots,N_p.
\]
With this notation, we formulate the computation of the electric field as a series of matrix-vector multiplications.
\begin{algorithm}
\caption{Grid-free evaluation of the electric field with quadratic complexity.}
\label{alg:electric-field-quadratic}
\begin{algorithmic}
\REQUIRE Mesh $\mathcal T=\left\{ \Gamma_k \right\}_{k=1}^{N_\Gamma}$, particles $(x_i)_{i=1}^{N_p}$ and matrices $V_h, K_h, M_h, W_h$.
\STATE $\underline N_0 \gets 1/\beta \, \Phi_{\mathcal T} \, w_q$
\IF{$\Gamma_N \neq \emptyset$}
\STATE $\underline N_1 \gets (K_h^\top - \frac{1}{2} M_h^\top) V_h^{-1} \underline N_0$.
\ENDIF
\IF{$\Gamma_N \neq \Gamma$}
\STATE Solve~\eqref{eq:BEM-block-system} for $\underline \phi_h$ and $\underline t_h$.
\ELSE
\STATE Solve~\eqref{eq:BIE_LGSNeumannDirichletVF} for $\underline \phi_h$.
\ENDIF
\FOR{$k=1$ \TO $3$}
\STATE $\underline E_k \gets -S_k \underline{t}_h + D_k \underline{\phi}_h - 1/\beta F_k w_q$
\COMMENT{Note the sign change due to $E = -\nabla \phi$}
\ENDFOR
\RETURN $\underline E_1, \underline E_2, \underline E_3$.
\end{algorithmic}
\end{algorithm}
The matrices in \cref{alg:electric-field-quadratic} are densely populated, see~\eqref{eq:galerkin-bem-matrices},~\eqref{eq:field-boundary-integrals}, and~\eqref{eq:field-free-space}.
Therefore, the algorithm scales like $\mathcal O(N_\Gamma^2 + N_\Gamma N_p + N_p^2)$ with a preprocessing step in the order of $\mathcal O(N_\Gamma^3)$ for computing
the Cholesky decomposition of the single layer operator $V_h$.

%% file: hierarchical.tex
\section{Hierarchical approximation}\label{sec:hierarchical}

In this section, we discuss how to find approximations to the fully populated matrices needed
for computation of the electric field in~\cref{alg:electric-field-quadratic}.

A direct evaluation is both quadratic in memory and computational time,
which can be large, even for a relatively small number of discretisation parameters.
With the special structure of most of the matrices, it is possible to reduce 
storage requirements and computational costs to linear complexity
by means of hierarchical approximations of dense matrices,
called $\mathcal{H}^2$--matrices.

The matrices given in \cref{sec:BEM-Galkerin} and \cref{sec:BEM-Vlasov-Poisson}
fit in a larger framework of integrals and point evaluations of a general kernel function.
To unify the treatment within the $\mathcal H^2$--technique, let us introduce a more general notation.
For the rest of this section we fix two index sets $\mathcal{I}$ and $\mathcal{J}$,
with associated sets $X \subset \mathbb{R}^3$ and $Y \subset \mathbb{R}^3$,
representing particles, nodes or triangles of the surface mesh.

The matrices $A \in \mathbb{R}^{\mathcal{I} \times \mathcal{J}}$ 
arising from \eqref{eq:field-free-space}
are point evaluations of a kernel function $k$, so-called Nystr\"om matrices,
\begin{equation}\label{eq:collocation-matrix}
A[i,j] = k(y_j, x_i), \quad i \in \mathcal{I}, \, j \in \mathcal{J},
\end{equation}
where $k$ is the fundamental solution~\eqref{eq:fundamental-solution}
or one component of its gradient, and
\linebreak[3] $X = (x_i)_{i \in \mathcal{I}}$, $Y = (y_j)_{j \in \mathcal{J}}$.
The Galerkin-type BEM matrices from~\eqref{eq:galerkin-bem-matrices}
have the form
\begin{equation}\label{eq:bem-matrix}
A[i,j] = \int_\Gamma \int_\Gamma
k(x, y) \varphi_j(y) \psi_i(x) \, \text ds_y \, \text ds_x,
\quad i \in \mathcal{I}, \, j \in \mathcal{J},
\end{equation}
with trial functions $(\varphi_j)_{j \in \mathcal{J}}$,
whose supports are in $Y \subset \Gamma$
and test functions
$(\psi_i)_{i \in \mathcal{I}}$ with supports in $X \subset \Gamma$.
Again, $k$ denotes the fundamental solution~\eqref{eq:fundamental-solution}
or its normal derivative.

The $\mathcal H^2$-matrix approximation~\cite{Boerm2010,Hackbusch2015}
is a tree-based data structure which exploits low-rank factorisations of matrix blocks,
\begin{equation}\label{eq:low-rank-h2}
A\left.\right|_{\sigma \times \tau} \approx V \Sigma U^\top,
\end{equation}
where $\sigma \subset \mathcal I$, $\tau \subset \mathcal J$ represent
parts of $X$ and $Y$,respectively, which are far apart, a term specified in \cref{def:admissibility}.
Furthermore, $V \in \mathbb R^{\sigma \times r}$,
$U \in \mathbb R^{\tau \times r}$, $\Sigma \in \mathbb R^{r \times r}$.
The dimension of $\Sigma$, $r \in \mathbb N$, is called the rank of the approximation.
To significantly reduce the storage requirements and the computational complexity,
$r \ll \max\{\#\sigma, \#\tau\}$ must hold.

The triple $(V, \Sigma, U)$ may be computed by a truncated Singular Value Decomposition.
Although this leads to the best compression rates, i.e. minimal storage complexity,
this method still has overall quadratic computational complexity.
It is therefore key to use a method which reduces both computational and storage complexity,
preferably to linear cost. Over the years, several methods have been proposed, most notable
Taylor expansion~\cite{Hackbusch2015}, multipole expansion~\cite{Carrier1988,Cheng1999,Greengard1987,Greengard1988,Greengard1997},
Adaptive Cross Approximation~\cite{Bebendorf2000, Bebendorf2008, BebendorfRjasanow2003, RjasanowSteinbach2007},
interpolation~\cite{BoermGrasedyck2004,Boerm2005a}, Hybrid Cross Approximation~\cite{Boerm2005b},
or the Green hybrid method~\cite{Boerm2016}.
For our study, we choose an approximation based on interpolation, although other methods would also work.
From a theoretical point of view, the complexity reduction from quadratic to linear can be best understood when utilising interpolation.
Practically, the interpolation scheme is readily implemented and very flexible as it works with point evaluations of the kernel function only.
Before we give details on the approximation scheme, we define an admissibility condition for the far field and describe the necessary
tree structures for the $\mathcal H^2$-format.
\begin{definition}
\label{def:admissibility}
Suppose $\sigma \subset \mathcal I$ and $\tau \subset \mathcal J$ with corresponding subsets $X_\sigma \subset X$, $Y_\tau \subset Y$.
The sets $\sigma$ and $\tau$ are are $\eta$--admissible
for $\eta > 0$ if
\[
\max\{\operatorname{diam}(X_\sigma),\,\operatorname{diam}(Y_\tau)\} 
\leq \eta \operatorname{dist}(X_\sigma, Y_\tau).
\]
\end{definition}

Searching for the optimal partition of $\mathcal{I} \times \mathcal{J}$
in a sense that equation \eqref{eq:low-rank-h2} holds for most
blocks with minimal rank $r$ is prohibitively expensive.
Therefore, the partition of $\mathcal{I} \times \mathcal{J}$, called
block cluster tree, is constructed via partitions of $\mathcal{I}$
and $\mathcal{J}$. These are given as cluster trees, see~\ref{alg:cluster-tree} for their construction.
They start with the whole index set in their root. Step by step,
further levels are added, where each level represents a disjoint union of the original set.
The splitting of a node of the cluster tree stops if the number of indices is below a given threshold $n_\text{min}$.
For computational ease and also for the interpolation points used later,
axis-parallel bounding boxes $B$ are associated to all nodes of the cluster tree.
This simplifies the computation of the admissibility condition a lot.
\begin{algorithm}
\caption{Construction of a cluster tree}\label{alg:cluster-tree}
\begin{algorithmic}
\REQUIRE Index set $\mathcal I$, geometry $X$, axis-parallel box $B \supset X$
\IF{$\# \mathcal I < n_\text{min}$}
\RETURN $T(\mathcal I) = (\mathcal I, B)$
\ELSE
\STATE Split $B$ into disjoint bounding boxes $B_1,\dots,B_{m}$
\FOR{$k=1$ \TO $n$}
\STATE Find maximal $\sigma_k \subset \mathcal I$ with $X_{\sigma_k} \subset B_k$
\STATE Add son to tree by constructing a cluster tree for $\sigma_k$, $X_{\sigma_k}$, $B_{\sigma_k} = B_k$.
\ENDFOR
\ENDIF
\RETURN Cluster tree $T(\mathcal I)$.
\end{algorithmic}
\end{algorithm}
The boxes in \cref{alg:cluster-tree} are usually split via principal component analysis or cardinality splitting.
If one splits the boxes in each direction, the well-known octree is recovered.

The block cluster tree contains the partition of $\mathcal I \times \mathcal J$ into non-admissible (near field) and admissible (far field) blocks.
It is defined as the cluster tree of $\mathcal I \times J$ with respect to \cref{def:admissibility}.
Its construction is given in \cref{alg:block-cluster-tree}.
\begin{algorithm}
\caption{Construction of a block cluster tree}
\label{alg:block-cluster-tree}
\begin{algorithmic}
\REQUIRE cluster trees $T(\mathcal I)$, $T(\mathcal J)$ for geometries $X, Y$.
\STATE $b \gets \operatorname{isAdmissible(\mathcal I, \mathcal J)}$ \COMMENT{Use condition from \cref{def:admissibility}.}
\IF{$b$ \OR $\operatorname{sons}(T(\mathcal I)) = \emptyset$ \OR $\operatorname{sons}(T(\mathcal J)) = \emptyset$}
\RETURN $T(\mathcal I, \mathcal J) = \emptyset$
\ELSE
\FOR{$\sigma$ \textbf{in} $\operatorname{sons}(T(\mathcal I))$}
\FOR{$\tau$ \textbf{in} $\operatorname{sons}(T(\mathcal J))$}
\STATE Add son by constructing a block cluster tree for $\sigma, \tau, X_\sigma, Y_\tau$.
\ENDFOR
\ENDFOR
\ENDIF
\RETURN $T(\mathcal I \times \mathcal J)$
\end{algorithmic}
\end{algorithm}

With the characterisation of the admissible blocks, we are able to give
the algorithm of low-rank approximations~\eqref{eq:low-rank-h2}
by means of interpolation.
Suppose $\sigma \times \tau \subset \mathcal I \times \mathcal J$
to be an admissible block.
For $r = m^3 \in \mathbb N$, $(z_k)_{k=1}^r \subset [-1,1]^3$ denotes
tensorised one-dimensional Chebyshev nodes,
\[
    \cos\left( \frac{2\ell-1}{2m}\pi \right), \quad \ell=1,\dots,m.
\]
These reference nodes are mapped to the boxes $B_\sigma$, $B_\tau$,
defining nodes $x^{(\sigma)}$, $y^{(\tau)}$, respectively,
which then are used to define the tensorised Lagrange polynomials
$(L^{\sigma}_k)_{k=1}^r$ and $(L^{\tau}_{k})_{k=1}^r$.
We now expand the kernel from equations~\eqref{eq:collocation-matrix}
or~\eqref{eq:bem-matrix} into the Lagrange basis,
\[
k(y,x) \approx \sum\limits_{k=1}^r \sum\limits_{\ell=1}^r
L^{(\sigma)}_{k}(x) 
k(y^{(\tau)}_{\ell}, x^{(\sigma)}_k)
L^{(\tau)}_\ell(y).
\]
Plugging this ansatz into~\eqref{eq:collocation-matrix}, we obtain
\[
A\left.\right|_{\sigma \times \tau} \approx V_\sigma \Sigma_{\sigma \times \tau} U_\tau^\top,
\]
where
\[
V_\sigma[i, k] = L^{(\sigma)}_k(x_i), \quad U_\tau[j, \ell] = L^{(\tau)}_\ell(y_j),
\]
and
\[
\Sigma_{\sigma \times \tau}[k, \ell] = k(y^{(\tau)}_{\ell}, x^{(\sigma)}_k),
\]
for $i \in \sigma, j \in \tau$, $k,\ell = 1,\dots,r$.
Note that the interpolation matrices $V_{\sigma}$ and $W_\tau$ only depend on the
corresponding cluster but not on the block cluster $\sigma \times \tau$.
When employing the interpolation-based approximation to Galkerin matrices~\eqref{eq:bem-matrix},
only the definitions of $V_\sigma$ and $U_\tau$ change. In that case,
\[
V_\sigma[i, k] = \int_\Gamma L^{(\sigma)}_k(x) \, \psi_i(x) \, \text ds_x, \quad 
U_\tau[j, \ell] = \int_\Gamma L^{(\tau)}_\ell(y) \, \varphi_j(y) \, \text ds_y,
\]
for $i \in \sigma, j \in \tau$, $k,\ell = 1,\dots,r$.
The families $(V_\sigma)_{\sigma \in T(\mathcal I)}$
and $(U_\tau)_{\tau \in T(\mathcal J)}$ are called cluster bases
for $T(\mathcal I)$ and $T(\mathcal J)$, respectively.
An important property of the cluster bases is that they are nested
in the following sense.
Assume $\sigma \in T(\mathcal I)$ to be a non-leaf node with sons
$\sigma', \sigma'',\dots$.

We begin with the observation that $L^{(\sigma)}\left.\right|_{X_{\sigma'}}$
and $L^{(\sigma')}$ both span the same polynomial space on $X_{\sigma'} \subset X_\sigma$.
Written out, we have
\[
L^{(\sigma)}_\ell(x) = \sum_{k=1}^r L^{(\sigma)}_\ell(x^{(\sigma')}_k) L^{(\sigma')}_k(x), \quad x \in X_{\sigma'},
\]
for $\ell=1,\dots,r$.
Therefore,
\[
V_\sigma\left.\right|_{\sigma' \times r} = V_{\sigma'} E_{\sigma', \sigma},
\]
where $E_{\sigma', \sigma}$ is called transfer matrix. Its entries are given by
\[
E_{\sigma,\sigma'}[k, \ell] = L^{(\sigma)}_\ell(x^{(\sigma')}_k), \quad k,\ell=1,\dots,r.
\]
Up to a permutation of the indices in $\sigma$, $V_\sigma$ can be written block-wise as
\[
V_\sigma = \begin{pmatrix}
V_{\sigma'} E_{\sigma, \sigma'} \\
V_{\sigma''} E_{\sigma, \sigma''} \\
\vdots
\end{pmatrix}.
\]
This special format reduces the required storage as the cluster basis
only depends on small transfer matrices. The full knowledge of $I$ and $X$ is only needed in the leafs of the cluster tree.
The nested structure of the cluster basis enables us to formulate a very efficient algorithm for matrix-vector multiplication.
It is split into three parts. In a first step, called forward transform, we iterate through $T(\mathcal J)$ and recursively collect the contributions of the transfer matrices.
Multiplication with the matrix entries takes place in the interaction phase. In a last step a apply the backward transform for $T(\mathcal I)$, which is similar to the first step.
Before we give the algorithms for these three parts, let us introduce an auxiliary vector. For $u \in \mathbb R^{\mathcal I}$ we define
\[
\hat u = \{u_\sigma \in \mathbb R^r: \sigma \in  T(\mathcal I)\}.
\]
With a precomputed cluster basis, the first step of the $\mathcal H^2$--matrix-vector multiplication in \cref{alg:fast-mvm} is the forward transform from \cref{alg:forward-transform}.
Followed by the interaction phase which couples the input vector and the output vector via point evaluations of the kernel function.
The result is then obtained by a backward transform applied to the output vector, see \cref{alg:backward-transform}.
The contribution of the fully assembled non-admissible blocks, i.e. the near field, is added in a final step.
Due to the nested structure of the cluster bases, the cluster depth in order of $\log \# \mathcal I$ or $\log \# \mathcal J$
disappears from the complexity estimate. The proven storage and computational complexity estimates read
\begin{lemma}\label{lem:h2-complexity}
Under mild assumptions on $T(\mathcal I)$ and $T(\mathcal J)$, the storage and therefore the computational complexity for
\cref{alg:fast-mvm} is
\[
\mathcal O(r(\# I + \# J)).
\]
\begin{proof}
For the proof and details on the assumptions on the cluster trees we refer the reader to~\cite{Hackbusch2015}.
\end{proof}
\end{lemma}
\begin{algorithm}
\caption{Forward transform}\label{alg:forward-transform}
\begin{algorithmic}
\REQUIRE cluster basis $(U_\tau)_{\tau \in T(\mathcal J)}$, vector $w \in \mathbb R^{\mathcal J}$, index set $\mathcal J$.
\IF{$\mathcal I$ is a leaf node}
\STATE $\hat w_{\mathcal J} \gets U_{\mathcal J}^\top w\left.\right|_{\mathcal J}$
\ELSE
\FOR{$\tau$ \textbf{in} $\operatorname{sons}(\mathcal J)$}
\STATE forward transform $(U_\tau)_{\tau \in T(\mathcal J)}$, $w$, $\tau$
\STATE $\hat w_{\mathcal J} \gets \hat w_{\mathcal J} + E_{\mathcal J, \tau}^\top \hat w_\tau$
\ENDFOR
\ENDIF
\RETURN $\hat w$
\end{algorithmic}
\end{algorithm}

\begin{algorithm}
\caption{Backward transform}\label{alg:backward-transform}
\begin{algorithmic}
\REQUIRE cluster basis $(V_\sigma)_{\sigma \in T(\mathcal I)}$, set of vectors $\hat u$, index set $\mathcal I$.
\IF{$\mathcal I$ is a leaf node}
\STATE $u\left|\right._{\mathcal I} \gets V_{\mathcal I} \hat u_{\mathcal I}$
\ELSE
\FOR{$\sigma$ \textbf{in} $\operatorname{sons}(\mathcal I)$}
\STATE $u\left|\right._{\mathcal I} \gets u\left|\right._{\mathcal I} + E_{\mathcal I, \sigma} u\left|\right._{\tau}$
\STATE backward transform $(V_\sigma)_{\sigma \in T(\mathcal I)}$, $\hat u$, $\sigma$
\ENDFOR
\ENDIF
\RETURN $u$
\end{algorithmic}
\end{algorithm}

\begin{algorithm}
\caption{$\mathcal H^2$--matrix-vector multiplication $u \gets u + \alpha A w$}\label{alg:fast-mvm}
\begin{algorithmic}
\REQUIRE $\alpha \in \mathbb R$, $A \in \mathbb R^{\mathcal I \times \mathcal J}$, $w \in \mathbb R^{\mathcal J}$, $u \in \mathbb R^{\mathcal I}$,
$(V_\sigma)_{\sigma \in T(\mathcal I)}, (U_\tau)_{\tau \in T(\mathcal J)}$.
\STATE Compute $\hat w$ by a forward transform from $(U_\tau)_{\tau \in T(\mathcal J)}$, $w$, $\mathcal J$.
\FOR{$\sigma$ \textbf{in} $T(\mathcal I)$}
\STATE $\hat u_\sigma \gets 0$
\ENDFOR
\FOR{$\sigma \times \tau$ \textbf{in} $T(\mathcal I \times \mathcal J)$ admissible}
\STATE $\hat u_\sigma \gets \hat u_\sigma + \alpha \, \Sigma_{\sigma \times \tau} \, \hat w_\tau$
\ENDFOR
\STATE $u \gets 0$
\STATE Compute $u$ by a backward transform from $(V_\sigma)_{\sigma \in T(\mathcal I)}$, $\hat u$, $\mathcal I$
\FOR{$\sigma \times \tau$ \textbf{in} $T(\mathcal I \times \mathcal J)$ \emph{not} admissible}
\STATE $u\left.\right|_{\sigma} \gets u\left.\right|_{\sigma} + \alpha A\left.\right|_{\sigma \times \tau} \, w\left.\right|_{\tau}$
\ENDFOR
\end{algorithmic}
\end{algorithm}

%% file: implementation.tex
\section{Notes on the implementation}\label{sec:implementation}

In this section we discuss our scheme with regard to its implementation
utilising hierarchical matrices. We also give an overview of the employed
software packages. All approximations with $\mathcal{H}^2$--matrices
in this section use polynomial interpolation as discussed in \cref{sec:hierarchical}.

At each time step, the system of boundary integral equations 
\eqref{eq:BEM-System} is solved to obtain
the Dirichlet and Neumann traces for the representation formula
\eqref{eq:BEM-representation-formula}.
The matrices from the discrete formulation \eqref{eq:galerkin-bem-matrices}
only depend on the discretisation of the boundary~$\Gamma$ but not
on the Dirichlet or Neumann boundary data or the positions of the particles.
Therefore they are computed in a preprocessing step and stored.
Afterwards, they are used for simulations with same geometry 
but possibly different boundary data or particle distributions.
All BEM matrices are approximated by $\mathcal{H}^2$--matrices.
For moderately sized problems,
we also compute the inverse of the single layer potential 
and approximate it by a $\mathcal{H}^2$--matrix.
Although this leads to cubic complexity in the number of triangles,
solving the linear system directly is faster than using an iterative method.
For larger problems, this is not feasible anymore.
We then apply a preconditioned conjugate gradient method,
see \cite{Steinbach2007} for preconditioning techniques
in case of BEM matrices.

The computation of the discrete Dirichlet trace of the Newton potential
$\underline N_0$ in \eqref{eq:newton-dirichlet-mvm} requires the $L^2$--projection onto the space of 
piecewise constant trial functions.
This can be formulated as a matrix-vector multiplication.
The computationally expensive matrix is efficiently 
approximated by a $\mathcal{H}^2$--matrix.
Let us fix quadrature rules for all triangles of the surface mesh.
Ideally all nodes lie on the edges of the triangles,
therefore reducing the number of function evaluation 
as triangles sharing a common edge also share quadrature nodes
and only differ in the weights.
Typical choices are the midpoints of the edges of the triangles
or the vertices of the triangles.
We collect all quadrature nodes in a global set $(y_j)_{j\in\mathcal{Q}}$
and denote the positions of the particles by $(x_i)_{i=1}^{N_p}$.
We can write
\begin{equation}\label{eq:dirichlet-newton-quadrature}
\underline N_0 \approx \frac{1}{\beta} M_{\mathcal Q} \Phi_{\mathcal Q} w_q,
\end{equation}
where $w_q = (q_0 |\Omega|/N_p)_{i=1}^{N_p}$ is the vector of weighted charges and
\[
\Phi_{\mathcal Q}[j, i] = U_\varepsilon(y_j, x_i) = \frac{1}{4 \pi} \frac{1}{| y_j - x_i |},
\quad j \in \mathcal{Q},~i=1,\dots,N_p,
\]
are the evaluations of the fundamental solutions of the particles
at the quadrature nodes and $M_{\mathcal Q}$ is a sparse matrix mapping the
global nodes $(y_j)_{j \in \mathcal{Q}}$ to the corresponding triangles,
multiplied with the quadrature weights.

The matrix $\Phi_{\mathcal Q}$ is a Nystr\"om matrix whose approximation
by $\mathcal{H}^2$--matrices is discussed in \cref{sec:hierarchical}.
The evaluation of $\underline N_0$ is reduced to linear complexity
in both the number of particles and the number of triangles.
The Neumann trace $\underline N_1$ is now readily computed
by the relation~\eqref{eq:neumann-newton-potential}
in linear complexity given the matrices are approximated by
$\mathcal{H}^2$--matrices.

For the computation of the electric field the gradient of the
representation formula is evaluated at the positions of the particles.
In order to efficiently apply the gradient, each evaluation of a component is
reformulated as a matrix-vector product, see \cref{alg:electric-field-quadratic}.
Due to the special structure of the fundamental solution~\eqref{eq:fundamental-solution}, these products are computed simultaneously.
Since
\[
-\nabla_x U(x,y) = \frac{x - y}{4 \pi |x-y|^3}, \quad x \neq y \in \mathbb R^3,
\]
the entries of the matrices $(S_k)_{k=1}^3$, $(D_k)_{k=1}^3$ and $(F_k)_{k=1}^{3}$ from equations~\eqref{eq:field-boundary-integrals} and~\eqref{eq:field-free-space}
have a common denominator.
Additionally, its computation is the most expensive part of the algorithm.
Therefore we compute the third power of the distance only once
and use this result for the computation of all matrices.
Again, these matrices fit into the general framework from \cref{sec:hierarchical}
and are efficiently approximated by $\mathcal{H}^2$--matrices reducing the
complexity of the matrix-vector multiplication from quadratic to linear
with respect to the number of particles.
Summarising, the quadratic \cref{alg:electric-field-quadratic} is transformed to an algorithm with
linear complexity by replacing all full matrices by their $\mathcal H^2$--approximations and
by performing the matrix-vector multiplication as described in \cref{alg:fast-mvm}.
Transferred to a computer program, the subroutines for a dense matrix-vector multiplication
simply have to be changed to their $\mathcal H^2$--matrix equivalents.
In this sense, using $\mathcal H^2$--matrices accelerates the algorithm without changing
important properties like the exact evaluation of the Coulomb force in the near field or the highly
accurate evaluation of the gradient of the representation formula, see \cref{lem:BEM-pointwise-error}.

In our scheme, the aforementioned $\mathcal{H}^2$--matrices
involving the positions of the particles are never fully built.
Instead we exploit their hierarchical structure and compute
the matrix-vector products on the fly.
Iterating through the block cluster tree and accumulating the
contribution of the admissible leafs,
the computation of the full matrices $\Phi_{\mathcal Q}, (S_k)_{k=1}^3, (D_k)_{k=1}^3$, and $(F_k)_{k=1}^3$
are reduced to the computation of the small leaf matrices.
Only storage for these small matrices is allocated which are freed after
a matrix-vector multiplication with parts of the vector $w_q$.
The positions of the particles change after each time step.
It is therefore necessary to rebuild the cluster tree, block cluster trees
and the cluster basis.
Although with a formal complexity of
$\mathcal{O}(N_p \log N_p)$,
the computational time is negligible compared to the computation
of the BEM gradient, see the timings in \cref{sec:numerics}.

The computation of the electric field relies heavily on an
efficient implementation of the hierarchical matrix format and tree-based data
structures.
We developed our code based on the H2Lib\footnote{%
The source code and further information can be found at \url{http://h2lib.org/}.}.
Written in the programming language C,  all basic data structures
and higher level routines like matrix-vector and matrix-matrix multiplication
or factorisation algorithms are available, as well as
a BEM module for the Laplace equation in three dimensions,
which is used in the subsequent computations.

%% file: numerics.tex
\section{Numerical examples}\label{sec:numerics}

In this section we present several numerical examples.
We begin with benchmarking the evaluation of the electric field
and conclude with physically motivated examples that
demonstrate classical plasma phenomena.

\subsection{Verification of linear complexity}

We numerically validate the linear scaling of the computational time
for the evaluation of the electric field at the positions of the particles.
The computation is split into four parts:
\begin{enumerate}
\item Building the cluster basis in $\mathcal{O}(N_p \log N_p)$,
\item computation of $\underline{N}_0$ according to~\eqref{eq:dirichlet-newton-quadrature} in linear complexity,
\item computation of the particle-particle force, see~\eqref{eq:field-free-space}
in linear complexity, and
\item evaluating the gradient of the representation formula~\eqref{eq:BEM-representation-formula} in linear complexity.
\end{enumerate}
For our tests, we triangulate the surface of the unit ball in $\mathbb{R}^3$
and uniformly distribute negatively charged particles inside the domain.
Appropriate nondimensionalisation is irrelevant for this test,
so we set all masses, charges and weights to unity.
Homogeneous Dirichlet boundary conditions are chosen for the electric potential.
We use $m=5$ interpolation nodes at each spatial direction for the
$\mathcal{H}^2$--matrix approximation.
The minimal cluster leaf size $n_\text{min}$ is $2 m^3$
and the admissibility constant $\eta$ is $2$.

\begin{figure}
\centering
\begin{minipage}{0.45\linewidth}
\includegraphics[width=\textwidth, page=1]{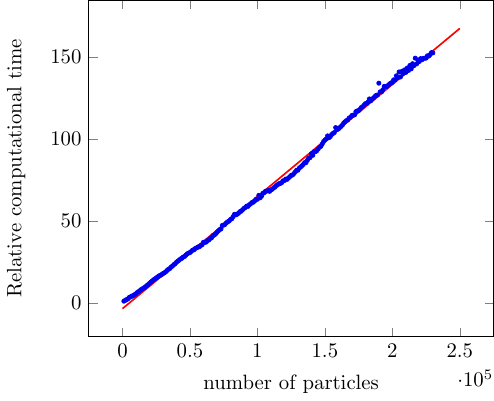}
\subcaption{Rebuilding cluster basis}\label{fig:timing_cluster}
\end{minipage}\hfill
\begin{minipage}{0.45\linewidth}
\includegraphics[width=\textwidth, page=2]{timing.pdf}
\subcaption{Newton potential}\label{fig:timing_newton}
\end{minipage}\\[0.5cm]

\begin{minipage}{0.45\linewidth}
\includegraphics[width=\textwidth, page=3]{timing.pdf}
\subcaption{field of the particles}\label{fig:timing_particle}
\end{minipage}\hfill
\begin{minipage}{0.45\linewidth}
\includegraphics[width=\textwidth, page=4]{timing.pdf}
\subcaption{representation formula}\label{fig:timing_repgrad}
\end{minipage}
\caption{
Single-core timings relative to 1\,000 particles
for building the cluster basis in \cref{fig:timing_cluster} ($2.7\cdot 10^{-3}$\,s), 
the evaluation of $\underline N_0$ in \cref{fig:timing_newton} ($2.3\cdot 10^{-3}$\,s),
the particle field in \cref{fig:timing_particle} ($5.3\cdot 10^{-3}$\,s)
and the gradient of the representation formula in \cref{fig:timing_repgrad} ($1.4\cdot 10^{-1}$\,s).
The number of triangles is 1\,280.
Computations were performend on an Intel Xeon Gold 6154@3\,GHz with \texttt{icc 19}. Relevant compiler flags
are \texttt{-Ofast -xHost}.}
\label{fig:relative-timings}
\end{figure}

\begin{figure}
\centering
\includegraphics[width=\linewidth]{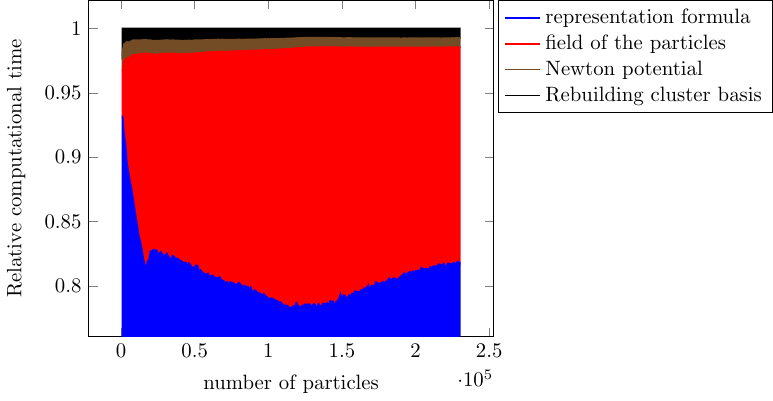}
\caption{Cumulative computational times
for 1\,280 triangles and varying number of particles.}
\label{fig:cumulative-timings}
\end{figure}

\begin{figure}
\centering
\includegraphics{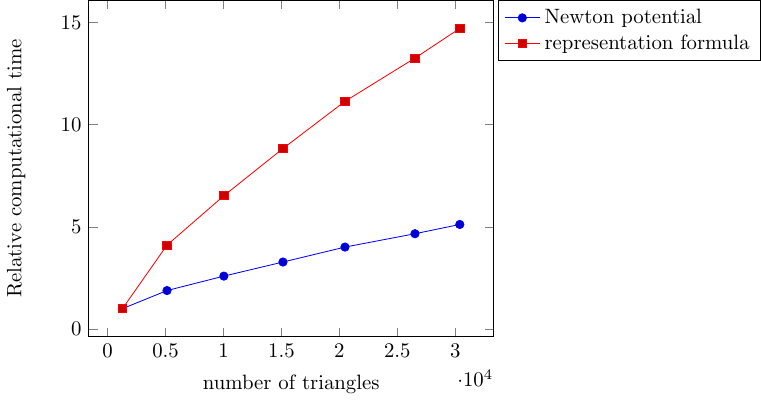}
\caption{Single-core timings for different number of triangles
relative to 1\,280 triangles for 10\,000 particles.
The base values are $1.9 \cdot 10^{-2}$\,s and $1.4$\,s for the Newton potential and representation formula, respectively.
For information on the CPU and the compiler, see the caption of \cref{fig:relative-timings}.}
\label{fig:bem-timings}
\end{figure}

\Cref{fig:relative-timings} shows the relative computational times
for a fixed mesh with varying number of particles.
The relative magnitudes of the different steps during the computation
of the electric field are given in \cref{fig:cumulative-timings}.
Although formally being of complexity $\mathcal{O}(N_p \log N_p)$,
we observe a linear scaling of the computation of the cluster basis.
Furthermore, the absolute timings are in the order of 100\,ms
making this part of the algorithm negligible compared to
rest of the algorithm which takes in the order of seconds.
The evaluation of $\underline{N}_0$ almost perfectly scales
linearly with the number of particles.
The evaluation of the gradient of the Newton potential and of the
representation formula follow a linear trend.
The constant hidden in the $\mathcal{O}$ notation
of \cref{lem:h2-complexity}
depends on the form of the block cluster tree.
As the particles are distributed randomly in the unit ball,
we cannot expect to obtain the same shape constant for the block
cluster tree for a large range of numbers of particles.
\Cref{fig:bem-timings} shows that the computation of 
the gradient of  the representation formula and of $\underline{N}_0$
scale linearly with the number of triangles.

\subsection{Physically motivated examples}

For most applications the plasma contains 
positively and negatively charged particles.
Usually, the positive charge exists of ionised atoms
and electrons form the negatively charged part.
Since the atoms are much heavier than the electrons
they are modelled as immobile. This gives rise
to a homogeneous positive background charge,
such that the system is electrically neutral from the outside.
The Poisson equation in~\eqref{eq:vlasov-poisson} changes
to
\[
-\Delta_x \phi = \frac{1}{\beta}\left[ 
1 - w \sum\limits_{i=1}^{N_p} \delta^\varepsilon_{x_i}
\right]
\]
with boundary conditions
\[
\begin{aligned}
\phi &= g_D & \text{on } \Gamma_D, \\
n_\Omega \cdot \nabla \phi &= g_N & \text{on } \Gamma_N.
\end{aligned}
\]
Note that the integral of the right-hand side over $\Omega$ is zero,
as $w = |\Omega| / N_p$.
A particular solution for the homogeneous background charge is
\[
\phi_b(x) = -\frac{1}{6\beta} | x |^2 \quad x \in \Omega.
\]
By subtracting traces of the particular solution $\phi_b$, we transform
the boundary value problem to
\[
\begin{aligned}
-\Delta_x \phi_e &= \frac{1}{\beta} w \sum\limits_{i=1}^{N_p} \delta_{x_i} \\
\phi_e &= g_D - g_b & \text{on } \Gamma_D, \\
n_\Omega \cdot \nabla \phi_e &= g_N - n_\Omega \cdot \nabla \phi_b & \text{on } \Gamma_N.
\end{aligned}
\]
The electric field is now obtained by
\[
E = -\nabla \phi_e - \nabla \phi_b.
\]
As $\phi_b$ is independent of the geometry and the distribution of the particles,
its evaluation and the evaluation of its gradient are grid-free, as well as
the computation of $\phi_e$.
Computations with background charge can be found in \cref{sec:plasma-oscillation}
and \cref{sec:plasma-sheath}.

\subsubsection{Accelerator}

\begin{figure}
\centering
\includegraphics{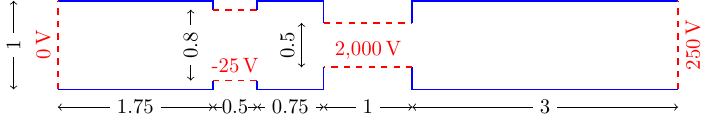}
\caption{Profile of the accelerator along the $x$-axis.
On the solid blue parts, homogeneous Neumann conditions are imposed.
The voltages along the dashed red lines indicate the value of the Dirichlet
boundary condition on these segments.}
\label{fig:accelerator-profile}
\end{figure}

\begin{figure}
\centering
\begin{minipage}{0.8\linewidth}

\begin{minipage}{\textwidth}
\centering
\includegraphics[width=\textwidth]{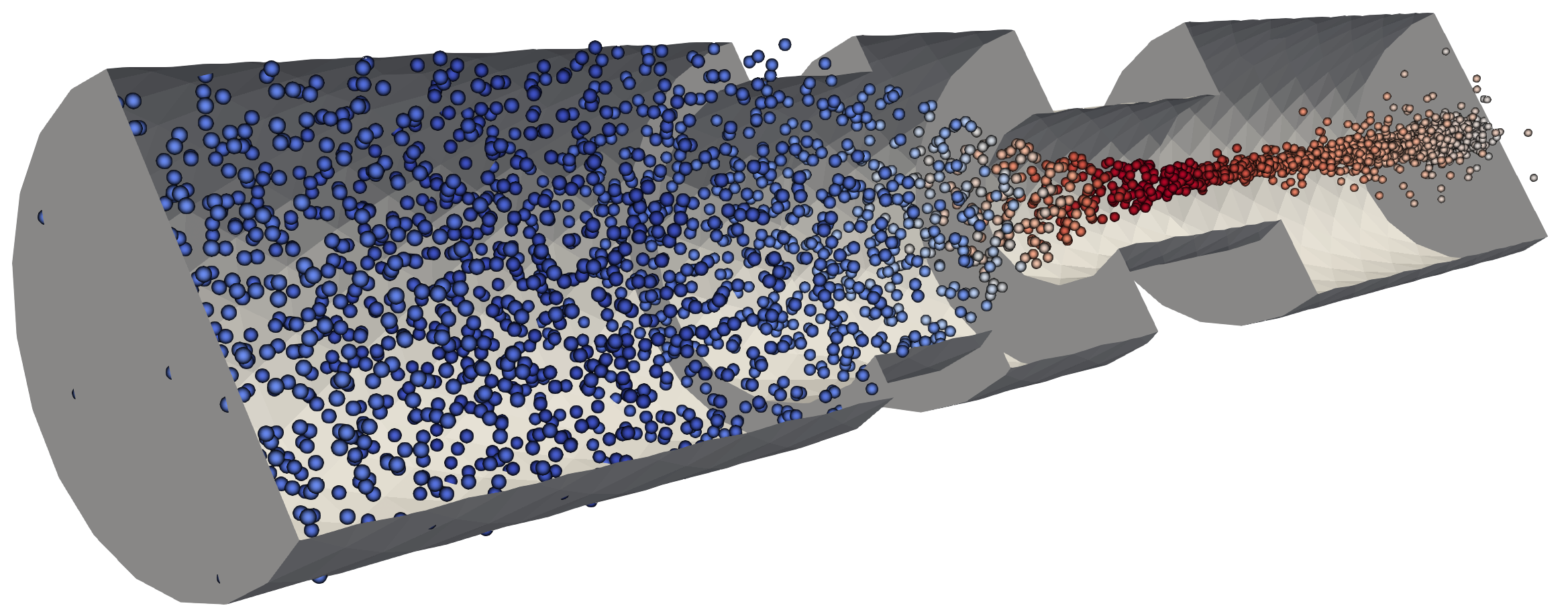}
\end{minipage}\\

\begin{minipage}{\textwidth}
\centering
\includegraphics[width=\textwidth]{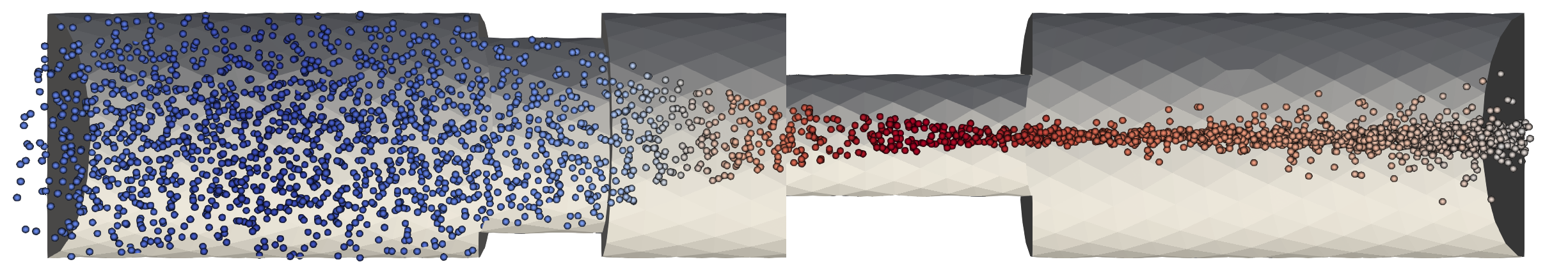}
\end{minipage}

\end{minipage}\hfill
\begin{minipage}{0.2\linewidth}
\centering
\includegraphics{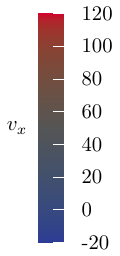}

\end{minipage}

\caption{Final distribution of 3\,000 particles inside the accelerator.
The colour indactes the velocity of the particles in $x$-direction.
}
\label{fig:accelerator-final}
\end{figure}

As a first example for non-trivial boundary conditions, we consider
an accelerator geometry, meshed with 8\,904 triangles.
The physically relevant parameters are $L_0 = 0.1 \, \text{m}$,
$n_0=10^{12} \, \text{m}^{-3}$ and $k_B T_0 = 1 \, \text{eV}$.
The profile of the rotationally symmetric accelerator
and the boundary conditions for the electric potential are depicted
in \cref{fig:accelerator-profile}.
Initially, 10\,000 particles are placed in the left cylinder with
a bulk velocity of 10 in positive $x$-direction and are absorbed at
the boundary.
Once they pass the first narrow, called screen, they are focused
such that they pass the second narrow, the accelerator, without being
absorbed by the boundaries.
The distribution of 3\,000 particles after 100 time steps with a time step
size of $10^{-3}$ is shown in \cref{fig:accelerator-final}.

\subsubsection{Plasma oscillations}\label{sec:plasma-oscillation}
\begin{figure}
\centering
\includegraphics{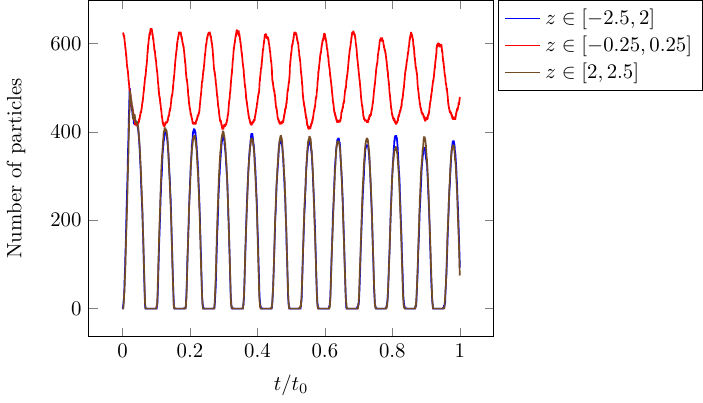}
\caption{Number of particles in three parts
of the cylinder over time for $n = 10 n_0$.}
\label{fig:particle-oscillation}
\end{figure}

\begin{figure}
\centering
\includegraphics{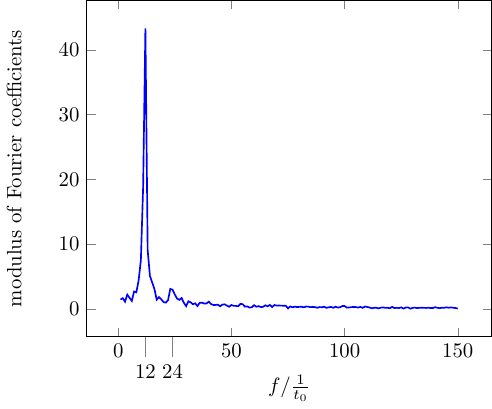}
\caption{Fourier spectrum of the number of particles in the middle of
the cylinder, the red line in \cref{fig:particle-oscillation}.
The constant mode is excluded from the spectrum.}
\label{fig:particle-oscillation-fourier}
\end{figure}

\begin{figure}
\centering
\includegraphics{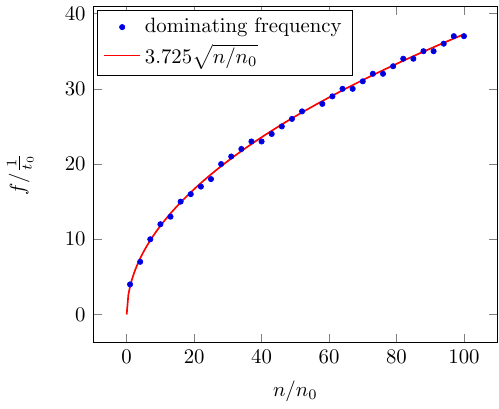}
\caption{Frequency of the oscillation of the number of particles
in the middle of the cylinder as a function of the electron
density.}
\label{fig:particle-oscillation-frequency}
\end{figure}

As a first example with a homogeneous background charge, we examine plasma
oscillations. 
The geometry is a cylinder along the $z$-axis with radius 1 and height 5,
centered in $0$. It is discretised with 2\,110 triangles.
The characteristic quantities are $L_0=0.1 \, \text{m}$,
$n_0 = 10^{12} \, \text{m}^{-3}$ and $k_B T_0 = 1 \, \text{eV}$.
5\,000~particles are distributed uniformly in a smaller cylinder of height 4 around
the centre of the geometry.
Their initial velocities are set to $0$.
The boundary is absorbing; at the bases we set homogeneous Dirichlet conditions
and homogeneous Neumann conditions on the rest.
Phyiscally, the latter boundary condition means that we impose a vanishing surface
charge density, in particular there is no net charge on this part of the boundary.
Mathematically, since the normal vectors point in radial direction,
the condition
\[
0 = \gamma_1 \phi = n_\Omega \cdot \gamma_0 \nabla \phi = - n_\Omega \cdot \gamma_0 E
\]
ensures that the field lines close to the boundary are parallel to the cylinder axis.
To prevent the particles from being absorbed at the lateral surface of the cylinder,
we add a constant magnetic field in the order of $10 \, \text{mT}$ along
the $z$-axis. The acceleration due to the magnetic field
is computed with the Boris scheme~\cite{Birdsall2004,Boris1970}
using a time step size of $10^{-4}$.
In an infinite system, the plasma oscillates with the plasma frequency
\[
\omega_p = \sqrt{\frac{n_0 e^2}{\varepsilon_0 m_e}},
\]
which depends only on the electron density.
As we simulate the plasma in a bounded domain, we cannot expect
the plasma to oscillate with the frequency $\omega_p$.
Instead, we validate that the frequency for the bounded domain is
still a function of the square root of $n_0$.
In order to do so, we vary the electron density $n$ from $n_0$ to $100 n_0$.
Counting the number of particles in three parts of the cylinder, 
$z \in [-2.5,2]$, $z \in [-0.25,0.25]$ and $z \in [2,2.5]$ at each time step,
we extract the dominating non-zero frequency after with the help
of the Discrete Fourier Transform.
The numbers of particles in the left, the middle and right part of the
cylinder for $n = 10 n_0$ is shown in \cref{fig:particle-oscillation}.
The distribution of the particles oscillates with dominating frequency
of $12$ in units of $1/t_0$, which corresponds to a angular frequency of
\[
\omega_c = 3.2\cdot10^{8} \, \frac{1}{\text{s}}
\]
in physical units. This is in the order of the plasma frequency
\[
\omega_p = \sqrt{\frac{10 n_0 e^2}{\varepsilon_0 m_e}} \approx
1.8\cdot10^{8} \, \frac{1}{\text{s}}.
\]
The spectra of the lines in \cref{fig:particle-oscillation}
only differ in magnitude, not in the positions of peaks.
Therefore, we only show the spectrum of the second line of
\cref{fig:particle-oscillation} in
\cref{fig:particle-oscillation-fourier}.
Repeating this several densities between $n_0$ and $100 n_0$ yields
\cref{fig:particle-oscillation-frequency},
from which the dependency of the frequency on the square root
of the density is clearly deduced.

\subsubsection{Plasma sheath}\label{sec:plasma-sheath}

\begin{figure}
\centering
\includegraphics{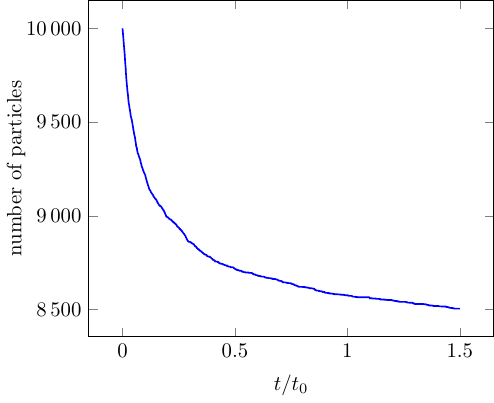}
\caption{Number of particles inside the sphere over time.}
\label{fig:sheath-number-particles}
\end{figure}

\begin{figure}
\centering
\includegraphics{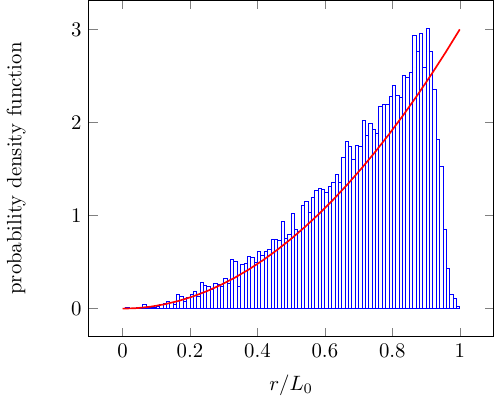}
\caption{Radial histogram of the final particle distribution inside the sphere.
The solid red line shows the probability density function of the uniform distribution.}
\label{fig:sheath-final-histogram}
\end{figure}

A classical nonlinear phenomenon in plasma physics is the formation of
sheaths, see the classical textbook~\cite{Chen2016}. For this example, we set
$L_0 = 0.1\,\text{m}$, $n_0 = 10^{13} \, \text{m}^{-3}$
and $k_B T_0 = 1 \, \text{eV}$.
We uniformly distribute 10\,000 particles
with velocity following a Maxwellian distribution with temperature~1
and bulk velocity~0 within the unit sphere,
which is discretised with 1\,280 triangles.
The particles are absorbed at the boundary; for the electric potential,
we impose homogeneous Dirichlet boundary conditions.
The system is evolved with a time step size of $10^{-3}$.
\cref{fig:sheath-number-particles} shows the number of
particles within the unit sphere as a function of time.
At the beginning, the fastest particles leave the sphere,
giving rise to a positive charge at the boundary.
With the growing potential barrier, the particles are excluded
from a thin area near the boundary, the so called sheath,
and are confined inside the sphere.
\Cref{fig:sheath-final-histogram} includes
the final radial distribution function of the particles inside the sphere
and the analytical radial distribution function of a uniformly distributed
random variate inside the unit sphere.
While the final positions are still uniformly distributed up to a radius
of approximately $0.6$, the distribution strongly deviates from the uniform
distribution especially close to radii of $1$, where
it suddenly drops to $0$.

\subsection{Summary}
To summarise, the numerical examples show that we are capable to simulate
important non-linear plasma phenomena like plasma oscillations or the
formation of sheaths.
The results also match available theoretical predictions.
Furthermore, the numerical study demonstrates the linear complexity of
our method and its applicability on three-dimensional domains
with mixed boundary values.
The efficiency and flexibility of our approach open the possibilities
for future simulations of complex problems in different plasma regimes.